\documentclass[3p]{scrartcl}

\usepackage{titling}
\usepackage[hmarginratio=1:1,top=32mm,left=15mm,columnsep=1cm]{geometry} %
\usepackage{array}
\usepackage{color}
\usepackage{tabularx}
\usepackage{graphicx}
\usepackage[svgnames,table,xcdraw]{xcolor}
\usepackage{amsmath}
\usepackage{amsthm}
\usepackage{amssymb}
\usepackage{amsfonts}
\usepackage{comment}
\usepackage{appendix}
\usepackage{bm}
\usepackage{multirow}
\usepackage{relsize} 
\usepackage{fourier} 
\usepackage{epstopdf}
\usepackage[font=small,labelfont=bf]{caption}	
\usepackage{multicol,booktabs}
\usepackage{nomencl}	
\usepackage{pbox} 
\usepackage{subcaption} 
\makenomenclature

\graphicspath{{../latex/draft/}}	

\newcommand\BibTeX{{\rmfamily B\kern-.05em \textsc{i\kern-.025em b}\kern-.08em
T\kern-.1667em\lower.7ex\hbox{E}\kern-.125emX}}

\usepackage{hyperref}
\hypersetup{
    colorlinks=true,
    linkcolor=Blue, 
    citecolor=DarkRed, 
    urlcolor=blue  } 

\newcommand{\pd}{\mathcal{\partial}}

\newcommand{\q}{\boldsymbol{q}}
\newcommand{\x}{\mbf{x}}
\newcommand{\y}{\mbf{y}}

\newcommand{\mbf}[1]{\mathbf{#1}}			%

\newcommand{\tref}{\text{ref}}
\newcommand{\tkin}{\text{kin}}				%

\newcommand{\Smix}{S}

\newcommand{\uphi}{{v_\theta}} 
\newcommand{\wphi}{{w_\theta}}
\newcommand{\ur}{{v_r}} 
\renewcommand{\wr}{{w_r}}

\newcommand{\RS}{\text{RS}}

\newcommand{\bdm}{\begin{displaymath}}
\newcommand{\edm}{\end{displaymath}}

\newcommand{\bea}{\begin{eqnarray} }
\newcommand{\eea}{\end{eqnarray} }

\renewcommand{\div}{{\nabla \cdot}}

\newcommand{\vv}{{\boldsymbol{v}}}

\newcommand{\ww}{{\bm{w}}}

\newcommand{\nn}{{\bm{n}}}

\newcommand{\Id}{{\boldsymbol{I}}}

\newcommand{\mf}{c}			
\newcommand{\vf}{\alpha}	
\newcommand{\indd}[1]{_{#1}}
\newcommand{\indu}[1]{^{#1}}
\newcommand{\rrat}[1]{\frac{\rho\indd{#1,\tref}}{\rho\indd{\tref}}}

\newcommand*\samethanks[1][\value{footnote}]{\footnotemark[#1]}

\newtheorem{theorem}{Theorem}[section]

\newtheorem{corollary}[theorem]{Corollary}

\newtheorem{assumption}[theorem]{Assumption}
\newtheorem{definition}[theorem]{Definition}

\usepackage{letltxmacro}

\LetLtxMacro\davidsincludegraphics\includegraphics

\makeatletter

\def\@starttocbutdonotshowit#1{%
	\begingroup
	\makeatletter
	\if@filesw
	\expandafter\newwrite\csname tf@#1\endcsname
	\immediate\openout \csname tf@#1\endcsname \jobname.#1\relax
	\fi
	\@nobreakfalse
	\endgroup}

\newcommand{\listoffigurenumbernames}{%
	\@starttocbutdonotshowit{lfn}%
}

\DeclareMathAlphabet{\mathsfbi}{OT1}{\sfdefault}{bx}{sl}
\DeclareMathVersion{sfletters}
\SetSymbolFont{letters}{sfletters}{OML}{ntxsfmi}{b}{it}

\makeatletter
\newcommand{\mathbfsbilow}[1]{%
	\text{\mathversion{sfletters}$\m@th#1$}%
}
\DeclareRobustCommand{\tensor}[1]{%
	\begingroup
	\ifcat\noexpand #1\relax
	\edef\greek@test{\detokenize{#1}}%
	\edef\greek@test{\expandafter\@cdr\greek@test\@nil}%
	\edef\greek@test{\expandafter\@car\greek@test\@nil}%
	\edef\x{\the\lccode\expandafter`\greek@test}%
	\edef\y{\number\expandafter`\greek@test}%
	\ifnum\x=\y\relax
	\mathbfsbilow{#1}%
	\else
	\mathsfbi{#1}%
	\fi
	\else
	\mathsfbi{#1}%
	\fi
	\endgroup
}
\makeatother

\newfont{\numerikEleven}{ecrm1000}
\newfont{\numerikTen}{cmss10}
\newfont{\numerikNine}{cmss9}
\newfont{\numerikEight}{cmss8}

\begin{document}

\title{
	 \textbf{An implicit-explicit solver \\
for a two-fluid single-temperature model}
}
%
%

\author{
M{\'a}ria Luk{\'a}{\v{c}}ov{\'a}-Medvid'ov{\'a},\thanksgap{-0.5ex}\thanks{Institut f\"ur Mathematik, Johannes-Gutenberg-Universit\"at Mainz, Staudingerweg 9, 55099 Mainz, Germany,
    (\href{mailto:lukacova@mathematik.uni-mainz.de}{lukacova@mathematik.uni-mainz.de})}
\quad
Ilya Peshkov\,\thanks{Department of Civil, Environmental and Mechanical Engineering, University of Trento, Via Mesiano 77, I-38123 Trento, Italy (\href{mailto:ilya.peshkov@unitn.it}{ilya.peshkov@unitn.it})
}\
\quad
and
\quad
Andrea Thomann\samethanks[1]
\textsuperscript{,}
\thanksgap{-1.25ex}
\thanks{Universit\'e de Strasbourg, CNRS, Inria, IRMA, F-67000 Strasbourg, France,
    (\href{mailto:andrea.thomann@inria.fr}{andrea.thomann@inria.fr})}
\textsuperscript{,$\ast$}
}
\thanksmarkseries{arabic}

\maketitle 

\date{\today}
{\let\thefootnote\relax\footnote{\hspace{0.25cm}$^* $Corresponding author}}

\paragraph{Abstract:} We present an implicit-explicit finite volume scheme for two-fluid single-temperature flow in all Mach number regimes which is based on a symmetric hyperbolic thermodynamically compatible description of the fluid flow.
The scheme is stable for large time steps controlled by the interface transport and is computational efficient due to a linear implicit character.
The latter is achieved by linearizing along constant reference states given by the asymptotic analysis of the single-temperature model.
Thus, the use of a stiffly accurate IMEX Runge Kutta time integration and the centered treatment of pressure based quantities provably guarantee the asymptotic preserving property of the scheme for weakly compressible Euler equations with variable volume fraction.
The properties of the first and second order scheme are validated by several numerical test cases.




\section{Introduction} \label{sec:introduction}
\setcounter{footnote}{0}

In continuum mixture theory, the constituents of a multiphase system, also called mixture, are present at every material element even if an element represents a pure phase.
This approach is applicable to model both situations -- the case of miscible \cite{PKG_Book2018} and immiscible \cite{Saurel2018,Chiocchetti2021} multicomponent systems.
The material interfaces, if present, are  zones of rapid but smooth changes
of a parameter distinguishing the phases of the mixture, typically the volume or mass fraction
.

Despite the fact that almost any application in science and engineering deals with multiphase systems and there is an obvious need for a consistent and reliable mathematical model to describe such multicomponent systems, the continuum mixture theory is far from being complete and no widely
accepted model exists.
Perhaps, the most widely used approach is based on equations for every individual constituent of the system, i.e. phase mass balance, phase momentum balance, phase energy balance, etc.
The key problem here is to find a closure for this system of equations which is represented by the coupling terms describing the exchange of mass, energy, and momenta between the mixture constituents.
Note that a first-principle theory to provide such a closure is currently not available. Consequently, various heuristic and phenomenological approaches are used.
The Baer-Nunziato (BN) model \cite{BaeNun1986} introduced in 1986 is a representative of this class of mathematical formulations and since then an active line of research has been done to adapt it to various applications, \cite{SauAbg1999,FavrGavr2012,Saurel2018} and recently to low Mach number flows \cite{ReAbg2022}.

In this paper, we deal with another class of governing equations for mixtures, here in the single-temperature simplification, which represents an attempt to build a mixture theory based on the first-principle reasoning.
The equations belong to the class of the so-called Symmetric Hyperbolic Thermodynamically Compatible (SHTC)
equations \cite{God1961,GodRom1995,GodRom2003}.
The key ingredients are the variational principle and the second law of thermodynamics.
The variational principle is used to deduce a reversible part of the evolution equations that is subject to entropy conservation.
The second law of thermodynamics yields an irreversible part of a model and controls entropy production.
In contrast, to the BN class of mixture models, the governing equations of the SHTC model are formulated not directly in terms of the phase quantities
but mainly in terms of mixture quantities such as mixture mass density, mixture momentum, mixture energy, etc.
The SHTC equations can be rewritten in terms of the phase balance equations, i.e. in a BN form, see \cite{RomResTor2007,RomBelPes2016}.
In this way, a new term appears in the phase balance equations, that is usually missing in the BN-type models, see \cite{RomBelPes2016}.
The latter can be identified as lift forces \cite{Drew1987} acting on a rotating fluid element of one phase immersed in another one.

In this work, we are interested in a \emph{single}-temperature SHTC mixture model \cite{Romenski2011} which is a special case of the two-velocity, two-pressure, two-entropy SHTC model of two-phase flows derived in \cite{RomResTor2007,RomDriTor2010}.
In \cite{ThoDum2023} the full model has been numerically solved based on an explicit time integration.
The difficulty lies in handling the two entropy balance laws of the full two-fluid SHTC model and only one energy conservation law.
One the other hand, applications such as sediment transport, granular flows or aerosol transport can be modeled with the single temperature approach resulting in one mixture entropy balance law associated to one energy conservation law.
These applications lie within weakly compressible flow regimes, characterized by small Mach numbers.
A severe difficulty in the construction of a numerical scheme applied to weakly compressible flow
regimes is posed by the scale differences between acoustic and material waves.
The focus of the numerical simulation usually lies on the evolution of the slower material waves following the two-fluid interface for which a time step controlled by the local flow speed is sufficient.
The time step of an explicit scheme, as proposed in \cite{RomDriTor2010,RomTor2004} for compressible two-phase flow in the SHTC framework, is bounded by the smallest Mach number.
This leads to very restrictive time steps in the low Mach number regime and consequently to long computational times, especially when long time periods are considered.
This problem can be overcome by considering implicit-explicit (IMEX) time integrators, where fast waves are treated implicitly leading to the Courant-Friedrichs-Lewy (CFL) condition that is restricted only by the local flow velocity.
It allows larger time steps while keeping the material waves well resolved.
Additionally, an implicit treatment of the associated stiff pressure terms, which trigger fast acoustic waves, has the advantage that centered finite differences can be applied without loss of stability while guaranteeing a Mach number independent numerical diffusion of the scheme, see e.g. \cite{Dellacherie2010,GuiVio1999,Klein1995} for a discussion on upwind schemes.
Indeed, the correct amount of numerical diffusion is crucial to obtain the so-called asymptotic preserving (AP) schemes \cite{Jin1999}.
Since the flow regime of the two-phase flow considered here is characterized by two potentially distinct phase Mach numbers, different singular Mach number limits can be obtained depending on the constitution of the mixture.
For their formal derivation we apply asymptotic expansions, as done for the (isentropic) Euler equations, see \cite{CorDegKum2012,DegTan2011,Dellacherie2010,GuiVio1999,KlaiMaj1981,Klein1995,ParMun2005}and the references therein.
We refer the reader to our recent work on the isentropic SHTC two-fluid model \cite{LukPupTho2022}.
To obtain physically admissible solutions, especially in the weakly compressible flow regime, the numerical scheme has to yield correct asymptotic behavior.
This means a uniformly stable and consistent approximation of the limit equations as the Mach numbers tend to zero.

The profound knowledge of the structure of well-prepared initial data can be used to construct an AP scheme by applying a reference solution (RS)-IMEX approach.
This approach was successfully applied to construct AP schemes for the (isentropic) Euler equations \cite{BisLukYel2017,KaiSchSchNoe2017,KucLukNoeSch2021,ZeiSchKaiBecLukNoe2020}
and isentropic two-fluid flow \cite{LukPupTho2022}.
This leads to a stiff linear part which is then treated implicitly whereas the nonlinear higher order terms are integrated explicitly respecting the asymptotes in the low Mach number limit.
By doing this, nonlinear implicit solvers can be avoided which are computationally costly.

The paper is structured as follows. In Section \ref{sec:model}, we briefly recall the model and give its non-dimensional formulation.
For well-prepared initial data, we analyze its singular Mach number limit towards the incompressible Euler equations with variable density.
Using the knowledge of the limit reference state, we construct first a semi discrete scheme in Section \ref{sec.SemiDiscrScheme} and derive a fully discrete scheme in Section \ref{sec.FullDiscrScheme}.
The construction of higher order schemes within this framework is shortly discussed, too.
Further, the AP property of the scheme is proven in Section \ref{sec:APanalysis}. Finally, in Section \ref{sec:NumRes}, a series of 1D and 2D test problems is presented to numerically verifying the convergence of the proposed scheme and its behavior in compressible and weakly compressible flow regime.

\section{Single temperature two-fluid flow} \label{sec:model}
In this section we recall the SHTC two-fluid model derived in
\cite{RomResTor2007,RomDriTor2010}.
We concentrate on the model in the thermal equilibrium regime \cite{Romenski2011} which is a legitimate approximation for many applications mentioned in the Introduction.
Thus, we deal with the mixture of two
fluids in which every material element (control volume) is characterized by the temperature $T$
with $T = T_1 = T_2$, where the lower indices denote the respective phase $l=1,2$. Moreover, we
assume that every material element of volume $\mathcal{V}$ and mass $\mathcal{M}$ is occupied by both fluids, i.e. $ \mathcal{V}
= \nu\indd{1} + \nu\indd{2}$ and $\mathcal{M} = m\indd{1} + m\indd{2}$, with $\nu\indd{l}$ and $m\indd{l}$
being the volume and mass of the $l$-th phase in the control volume $V$.
However, to characterize the fluid content in a control volume, it is convenient to use
non-dimensional scalars: the volume fractions $\vf\indd{l}$ and mass fractions $\mf\indd{l}$
defined as
\begin{equation}\label{eq:vf.mf}
    \vf\indd{l} = \frac{\nu\indd{l}}{\mathcal{V}},
    \qquad
    \mf\indd{l} = \frac{m\indd{l}}{\mathcal{M}} = \frac{\varrho\indd{l}}{\rho} =
    \frac{\vf\indd{l}\rho\indd{l}}{\rho},
\end{equation}
where
\begin{equation}\label{eq:rho}
    \rho = \frac{\mathcal{M}}{\mathcal{V}} = \varrho\indd{1} + \varrho\indd{2} = \vf\indd{1}
    \rho\indd{1} + \vf\indd{2}\rho\indd{2}
\end{equation}
is the mass density of the mixture,
$\varrho\indd{l}$ is the mass density of the $l$-th phase in the control volume $V$, and
$\rho\indd{l}$ is the mass density of the $l$-th phase.
The volume and mass fractions obey the constraints
\begin{equation}\label{eq:constr1}
    \vf\indd{1} + \vf\indd{2} = 1,
    \qquad
    \mf\indd{1} + \mf\indd{2} = 1.
\end{equation}
Moreover, each phase is equipped with its own velocity field $\vv\indd{l} \in \mathbb{R}^d$, where
$d$
denotes the space dimension, and the mixture control volume is assume to have the velocity defined
as the center of mass velocity, i.e. as the weighted average given by
\begin{equation}
    \vv = \mf\indd{1}\vv\indd{1} + \mf\indd{2}\vv\indd{2}.
\end{equation}
The mixture momentum $\rho \vv = \varrho\indd{1}\vv\indd{1} +
\varrho\indd{2}\vv\indd{2}$ is equal to the sum of the phase momenta.
Additionally, one needs to characterize the relative motion of the phases which, in the SHTC
theory, is done using the relative velocity field
\begin{equation}\label{key}
    \ww = \vv\indd{1} - \vv\indd{2}.
\end{equation}
For each phase, an entropy $s\indd{l}$ and internal energy $e\indd{l}(\rho\indd{l},s\indd{l})$ is
prescribed yielding the phase pressures
\begin{equation}
    p\indd{l} = \rho\indd{l}^2 \frac{\pd e\indd{l}}{\pd \rho\indd{l}}, \quad l = 1, 2.
\end{equation}
We consider an {\em ideal gas} equation of state (EOS) given in terms of the respective density and single temperature resulting in
\begin{equation}
    \label{eq.EOS}
    s\indd{l}(\rho\indd{l},T) = c\indd{v,l} \log\left( \frac{T}{T_{0,l}}\left(\frac{1}{\rho\indd{l}}\right)^{\gamma\indd{l} -1}\right), \quad T\indd{0,l} = \frac{1}{(\gamma\indd{l} - 1)c\indd{v,l}}, \quad e\indd{l}(\rho\indd{l}, T) = c\indd{v,l} T, \quad p\indd{l}(\rho\indd{l},T) = (\gamma\indd{l} - 1) c\indd{v,l} \rho\indd{l} T,
\end{equation}
where $\gamma\indd{l}$ denotes the ratio of specific heats and $c\indd{v,l}$ the specific heat at constant volume for each phase $l=1,2$.
To conclude the definition of the mixture state variables, we have the specific mixture internal energy $e =
\mf\indd{1}e\indd{1} + \mf\indd{2}e\indd{2}$, the mixture pressure $p = \rho^2 \frac{\pd e}{\pd
    \rho} = \vf\indd{1}p\indd{1} +
\vf\indd{2}p\indd{2}$ and the total energy density of the mixture given by
\begin{equation}
    \rho E = \rho e + \rho \frac{\|\vv\|^2}{2} + \rho\mf\indd{1}\mf\indd{2}\frac{\|\ww\|^2}{2}.
\end{equation}
The total mixture entropy reads $S = \mf\indd{1}s\indd{1} + \mf\indd{2}s\indd{2}$.
All state variables are summarized in the vector
\begin{equation}
    \label{eq.stateVar}
    \q = \left(\vf\indd{1}, \vf\indd{1}\rho\indd{1}, \vf\indd{2}\rho\indd{2}, \rho\vv, \ww, \rho E\right)^T.
\end{equation}
The SHTC model with a single temperature can be written in the following way
\begin{subequations}\label{eq.SHTC.div}
    \begin{eqnarray}
        && \frac{\pd \vf\indd{1}}{\pd t}+\vv \cdot \nabla \vf\indd{1} =-\frac{p\indd{1}-p\indd{2}}{\tau^{(\vf)}\rho},\label{eqn.alphaMS.div}\\[2mm]
        && \frac{\pd (\vf\indd{1}\rho\indd{1})}{\pd t}+\div (\vf\indd{1}\rho\indd{1} \vv\indd{1}) = 0,\label{eqn.contiMS1.div}\\[2mm]
        && \frac{\pd (\vf\indd{2}\rho\indd{2})}{\pd t}+\div (\vf\indd{2}\rho\indd{2} \vv\indd{2}) = 0,\label{eqn.contiMS2.div}\\[2mm]
        &&\frac{\pd (\rho \vv)}{\pd t}+\div
        (\rho \vv \otimes \vv + p \Id + \rho \mf\indd{1} \mf\indd{2} \ww \otimes  \ww )=0,
        \label{eqn.momentumMS.div}\\[2mm]
        &&\frac{\pd \ww}{\pd t}+\div\left(\left[ \ww \cdot \vv+\mu_1 - \mu_2 + (1-2\mf\indd{1}) \frac{\|\ww\|^2}{2}\right]\Id\right) + ( \nabla \times \ww ) \times \vv=-\dfrac{
            \mf\indd{1} \mf\indd{2} \ww }{\tau^{(w)}},\label{eqn.relvelMS.div}\\[2mm]
        &&
        \frac{\pd (\rho E)}{\pd t}+
        \div \left(\vv (\rho E + p)  +\rho ~\left[\ww \cdot \vv +
        \mu_1 - \mu_2 + (1-2\mf\indd{1}) \frac{\|\ww\|^2}{2} \right] \mf\indd{1} \mf\indd{2} \ww
        \right)=0.\label{eqn.energyMS.div}
    \end{eqnarray}
\end{subequations}
Here $\mu\indd{l} = e\indd{l} + \frac{p\indd{l}}{\rho\indd{l}} - s\indd{l} T, ~l=1,2$ denote the
chemical potentials.
In the above formulation, the volume fraction is advected in a non-conservative way with the fluid
flow $\vv$ balanced by a pressure relaxation source term.
The mixture mass is conserved  due to \eqref{eqn.contiMS1.div}, \eqref{eqn.contiMS2.div}, the momentum due to \eqref{eqn.momentumMS.div} and the total energy due to \eqref{eqn.energyMS.div}.
The relative velocity is not conserved and driven by the difference in the chemical potentials $\mu
= \mu\indd{1} - \mu\indd{2}$ and a friction source term.
The relaxation parameters $\tau^{{\vf}}$ and $\tau^{(\ww)}$ characterize the relaxation rates of the mixture towards pressure ($p\indd1=p\indd{2}$) and relative velocity
($\vv\indd{1}=\vv\indd{2}$) equilibrium.

Moreover, the model is equipped with the entropy balance law
\begin{equation}\label{eqn.entropyMS.div}
    \frac{\pd (\rho \Smix)}{\pd t}+\div (\rho \Smix \vv) = \Pi \geq 0
\end{equation}
with the entropy production term
\begin{equation}\label{eq:Pi}
    \Pi =
    \frac{1}{T \tau^{(\vf)}\rho^2}(p\indd{1}-p\indd{2})^2 + \frac{1}{T
        \tau^{(w)}\rho^2}\mf\indd{1}^2\mf\indd{2}^2
    \Vert\ww\Vert^2.
\end{equation}
For details on the derivation of the model and its thermodynamical properties we refer to
\cite{RomTor2004,RomResTor2007,RomDriTor2010,RomBelPes2016}.

Since each phase is equipped with a respective pressure and density, a sound speed for each phase $a\indd{l}$, as well as a mixture sound speed $a$ can be defined by
\begin{equation}
    \label{eq.soundspeeds}
    \left(a\indd{l}\right)^2 = \frac{\pd p\indd{l}}{\pd \rho\indd{l}}=\gamma\indd{l}\frac{p\indd{l}}{\rho\indd{l}} \quad \text{and} \quad
    a^2 = {\frac{\pd p}{\pd \rho}}=\mf\indd{1}\left(a\indd{1}\right)^2 + \mf\indd{2}\left(a\indd{2}\right)^2.
\end{equation}
Accordingly, a Mach number can be assigned to each phase.
As usual it is defined by the ratio between the flow velocity $\vv$ and the sound speed $a\indd{l}$.
In the case that the flow is characterized by (at least) one small Mach number, different scales arise in the model that yield stiffness in the governing equations \eqref{eq.SHTC.div}.
To obtain a better understanding of the scales which are present in the model, we will rewrite system \eqref{eq.SHTC.div} in a non-dimensional form.

\subsection{Non-dimensional formulation of the two-fluid model}

Let us denote the non-dimensional quantities by $(~\tilde \bullet~)$ and the corresponding
reference value by $\left(\bullet\right)\indd{\tref}$.
We assume that the convective scales are of the same order, i.e. $\vv\indd{l,\tref} = \vv_\tref = x_\tref/t_\tref$.
The ratio of the phase densities, however, can be large, especially when considering a mixture of a
light gas and liquid.
To take this potentially large difference into account, we define two different reference densities $\rho\indd{l,\tref}$.
Note, that the volume fraction and mass fractions are already non-dimensional quantities.
Further, we define two reference pressures $p\indd{l,\tref}$ from which we can compute the reference sound speeds $a_\tref$ and reference internal energies $e_\tref$ via the EOS \eqref{eq.EOS}.
They are given by
\begin{equation}
    \left(a\indd{l,\tref}\right)^2 = \frac{p\indd{l,\tref}}{\rho\indd{l,\tref}}, \quad e\indd{l,\tref} = \frac{p\indd{l,\tref}}{\rho\indd{l,\tref}}, \quad T\indd{\tref} = \frac{1}{\gamma\indd{l} (\gamma\indd{l} -1) c\indd{v,l}}\frac{p\indd{l,\tref}}{\rho\indd{l,\tref}}, \quad l = 1, 2.
\end{equation}
The dimensional state variables $\q$ are then expressed as the product of non-dimensional quantities and reference values as follows
\begin{equation}
    \label{eq.scaled.quantities}
    \rho\indd{l} = \tilde \rho\indd{l} \rho\indd{l,\tref}, \quad p\indd{l} = \tilde p\indd{l} p\indd{l,\tref}, \quad e\indd{l} = \tilde e\indd{l} \frac{p\indd{l,\tref}}{\rho\indd{l,\tref}}, \quad \mu\indd{l} = \tilde \mu\indd{l} \frac{p\indd{l,\tref}}{\rho\indd{l,\tref}}, \quad \vv\indd{l} = \tilde \vv\indd{l} \vv\indd{\tref}, \quad \vv = \tilde \vv \vv\indd{\tref}, \quad \ww = \tilde \ww \vv\indd{\tref}.\\
\end{equation}
Further, a respective reference Mach number $M\indd{l}$
is assigned to each phase $\quad l = 1, 2$
\begin{equation}
    \label{eq.Mach}
    M\indd{l} = \frac{\vv\indd{\tref}}{a\indd{l,\tref}}
    .
\end{equation}
Inserting expressions \eqref{eq.scaled.quantities} into the dimensional equations \eqref{eq.SHTC.div}, dropping the tilde $(~\tilde \bullet ~)$ and using \eqref{eq.Mach}, we obtain the following non-dimensional formulation

\begin{subequations}\label{eq.SHTC.div.nd}
    \begin{eqnarray}
        && \frac{\pd \vf\indd{1}}{\pd t}+\vv \cdot \nabla \vf\indd{1} =-\frac{1}{\tau^{(\vf)}\rho}\left(\rrat{1}\frac{p\indd{1}}{(M\indd{1})^2}-\rrat{2}\frac{p\indd{2}}{(M\indd{2})^2}\right),\label{eqn.alphaMS.div.nd}\\[2mm]
        && \frac{\pd (\vf\indd{1}\rho\indd{1})}{\pd t}+\div (\vf\indd{1}\rho\indd{1} \vv\indd{1}) = 0,\label{eqn.contiMS1.div.nd}\\[2mm]
        && \frac{\pd (\vf\indd{2}\rho\indd{2})}{\pd t}+\div (\vf\indd{2}\rho\indd{2} \vv\indd{2}) = 0,\label{eqn.contiMS2.div.nd}\\[2mm]
        &&\frac{\pd (\rho \vv)}{\pd t}+\div
        \left(\rho \vv \otimes \vv + \left(\vf\indd{1}\rrat{1}\frac{p\indd{1}}{(M\indd{1})^2}+
        \vf\indd{2}\rrat{2}\frac{p\indd{2}}{(M\indd{2})^2}\right)\Id + \rho \mf\indd{1}\mf\indd{2}
        \ww \otimes  \ww \right)=0,
        \label{eqn.momentumMS.div.nd}\\[2mm]
        &&\frac{\pd \ww}{\pd t}+\div\left(\left[ \ww \cdot \vv+\frac{\mu\indd{1}}{(M\indd{1})^2} - \frac{\mu\indd{2}}{(M\indd{2})^2} + \left(1-2\mf\right) \frac{\|\ww\|^2}{2}\right]\Id\right) + ( \nabla \times \ww ) \times \vv=-\dfrac{
            \mf\indd{1}\mf\indd{2}\ww }{\tau^{(w)}},\label{eqn.relvelMS.div.nd}\\[2mm]
        &&
        \frac{\pd (\rho E)}{\pd t}+
        \div \left(\vv \left(\rho E + \vf\indd{1}\rrat{1}\frac{p\indd{1}}{(M\indd{1})^2}+ \vf\indd{2}\rrat{2}\frac{p\indd{2}}{(M\indd{2})^2}\right)  \right) \\
        &&\hskip1.1cm+ \div \left(\rho ~\left[ \ww \cdot \vv+\frac{\mu\indd{1}}{(M\indd{1})^2} - \frac{\mu\indd{2}}{(M\indd{2})^2} + \left(1-2\mf\indd{1}\right) \frac{\|\ww\|^2}{2}\right] \mf\indd{1}\mf\indd{2}\ww \right)=0 \label{eqn.energyMS.div.nd}
    \end{eqnarray}
\end{subequations}
with the scaled total energy
\begin{equation}\label{eqn.E.nd}
    E = \mf\indd{1} \rrat{1}\frac{e\indd{1}(\rho\indd{1},T)}{(M\indd{1})^2} + \mf\indd{2} \rrat{2}\frac{e\indd{2}(\rho\indd{2},T)}{(M\indd{2})^2} +\frac{\Vert \vv \Vert^2}{2}+ \mf\indd{1}\mf\indd{2} \frac{\Vert \ww \Vert^2}{2}
\end{equation}
and the mixture density $\rho = \tilde \rho \rho\indd{\tref} = \vf\indd{1}\tilde\rho\indd{1}\rho\indd{1,\tref} + \vf\indd{2}\tilde \rho\indd{2}\rho\indd{2,\tref}$ with
$\tilde \rho = \vf\indd{1} \tilde \rho\indd{1}\rrat{1} + \vf\indd{2} \tilde \rho\indd{2} \rrat{2}.$
In the next section we introduce well-prepared initial data that will be used for the formal asymptotic analysis of \eqref{eq.SHTC.div.nd} in the low Mach number limits.

\subsection{Well-prepared data and low Mach number limits}
\label{sec.WellPrepLowMachLimits}
As we have seen from the Mach number definition \eqref{eq.Mach}, the difference in the flow regimes of the two phases depends mainly on the material constants $\gamma\indd{l}$ and $c\indd{v,l}$.
In particular, from the single temperature EOS \eqref{eq.EOS}, we obtain
\begin{equation}
    \frac{a\indd{1}^2}{\gamma\indd{1}(\gamma\indd{1}-1)c\indd{v,1}} = \frac{a\indd{2}^2}{\gamma\indd{2}(\gamma\indd{2}-1)c\indd{v,2}} \Leftrightarrow a\indd{1}^2 = \frac{\gamma\indd{1}(\gamma\indd{1}-1)c\indd{v,1}}{\gamma\indd{2}(\gamma\indd{2}-1)c\indd{v,2}} a\indd{2}^2
\end{equation}
and consequently, with $a\indd{l,\tref} = \sqrt{\gamma\indd{l} (\gamma\indd{l} - 1) c\indd{v,l} T\indd{\tref}}$, we find a direct relation between two Mach numbers
\begin{equation}
    \label{eq.Mach_C}
    M\indd{1} = \mathcal{C} M\indd{2}, \quad \mathcal{C} = \sqrt{\frac{\gamma\indd{2}(\gamma\indd{2}-1)c\indd{v,2}}{\gamma\indd{1}(\gamma\indd{1}-1)c\indd{v,1}}} > 0.
\end{equation}
In the following, for simplicity, we consider the case $M\indd{1} = M\indd{2} = M$, where $0 < M
\ll 1$, i.e. $\mathcal{C} = 1$.
The cases $\mathcal{C} > 1$ and $\mathcal{C} < 1$ can be treated in a similar manner.
For a full analysis of model \eqref{eq.SHTC.div.nd} in the isentropic case we refer the reader to \cite{LukPupTho2022}, where the singular limit for two different Mach numbers $1 \gg M\indd{1} > M\indd{2} > 0$ and $1 \approx M\indd{1} \gg M\indd{2} > 0$ are considered.

We proceed by expanding sufficiently smooth phase state variables with respect to $M$.
Note that the volume and mass fractions are non-dimensional quantities and are not expanded with respect to the Mach number.
\begin{align}\label{eq.Mach.expansion.phase}
    \begin{split}
        \rho\indd{l} &= \rho\indd{l,(0)} + M \rho\indd{l,(1)} + \mathcal{O}(M^2), \qquad l=1,2, \\
        T &= T\indd{(0)} + M T\indd{(1)} + M^2 T\indd{(2)} + \mathcal{O}(M^3),\\
        \vv &= \vv\indd{(0)} + M \vv\indd{(1)} + \mathcal{O}(M^2).\\
    \end{split}
\end{align}
Since the relative velocity is subject to a relaxation process, we set $\tau^{(\ww)} = M$ leading to the desired zero background relative velocity $\ww_{(0)} = 0$ in the limit, thus
\begin{equation}
    \label{eq.Mach.expansion.w}
    \ww = M \ww_{(1)} + \mathcal{O}(M^2).
\end{equation}

To obtain Mach number expansions also for the remaining thermodynamical quantities, we apply EOS \eqref{eq.EOS} which yields
\begin{align}
    \begin{split}
        p\indd{l} &= \left(c\indd{v,l}(\gamma\indd{l}-1)\rho\indd{l,(0)}T\indd{(0)}\right) + M c\indd{v,l}(\gamma\indd{l}-1)\left(\rho\indd{l,(0)}T\indd{(1)}+\rho\indd{l,(1)}T\indd{(0)}\right) + \mathcal{O}(M^2),\\
        \mu\indd{l} &= \left(T\indd{(0)}\left(c\indd{v,l}\gamma\indd{l} - s\indd{l}(\rho\indd{l,(0)},T\indd{(0)})\right)\right) + M \left(\frac{c\indd{v,l}(\gamma\indd{l}-1)\left(\rho\indd{l,(0)}T\indd{(1)}+\rho\indd{l,(1)}T\indd{(0)}\right)}{\rho\indd{l,(0)}}-T\indd{(1)}s\indd{l}(\rho\indd{l,(0)},T\indd{(0)})\right)  + \mathcal{O}(M^2), \\
        \rho\indd{l}e\indd{l} &= \left(c\indd{v,l}\rho\indd{l,(0)}T\indd{(0)}\right) + M c\indd{v,l}\left(\rho\indd{l,(0)}T\indd{(1)}+\rho\indd{l,(1)}T\indd{(0)}\right) + \mathcal{O}(M^2)
    \end{split}
\end{align}
and imply the following expansions
\begin{align}\label{eq.Mach.expansion.EOS}
    \begin{split}
        p\indd{l} &= p\indd{l,(0)} + M p\indd{l,(1)} + M^2p\indd{l,(2)} + \mathcal{O}(M^3), \\
        \mu\indd{l} &= \mu\indd{l,(0)} + M \mu\indd{l,(1)} + M^2\mu\indd{l,(2)} + \mathcal{O}(M^3), \\
        \rho\indd{l}e\indd{l} &= \left(\rho\indd{l}e\indd{l}\right)\indd{(0)} + M
        \left(\rho\indd{l}e\indd{l}\right)\indd{(1)} + M^2\left(\rho\indd{l}e\indd{l}\right)\indd{(2)}
        + \mathcal{O}(M^3). \\
    \end{split}
\end{align}
We insert the Mach number expansions \eqref{eq.Mach.expansion.phase},\eqref{eq.Mach.expansion.EOS} and $\rho_{(0)} = \vf^1 \rho^1\indd{(0)} + \vf^2 \rho^2\indd{(0)}$ in the non-dimensional formulation \eqref{eq.SHTC.div.nd} and sort by the equal order of the Mach number.
Terms of the order $\mathcal{O}(M^{-2})$ and $\mathcal{O}(M^{-1})$ arise in the relaxation source term of equation \eqref{eqn.contiMS1.div.nd} yielding
\begin{equation}\label{eqn.p.relax.cond}
    \rrat{2} p\indd{1,(0)} = \rrat{1} p\indd{2,(0)} \quad \text{and} \quad \rrat{2} p\indd{1,(1)} = \rrat{1} p\indd{2,(1)},
\end{equation}
as well as in the momentum equation \eqref{eqn.momentumMS.div.nd}
\begin{eqnarray}
    &&\nabla \left(\vf\indd{1} \rrat{1}p\indd{1,(0)} + \vf\indd{2}\rrat{2} p\indd{2,(0)}\right) = 0
    \Leftrightarrow \nabla p\indd{1,(0)} = 0, ~\nabla p\indd{2,(0)} = 0 \quad \text{and} \\
    &&\nabla \left(\vf\indd{1} \rrat{1}p\indd{1,(1)} + \vf\indd{2}\rrat{2} p\indd{2,(1)}\right) = 0
    \Leftrightarrow \nabla p\indd{1,(1)} = 0, ~\nabla p\indd{2,(1)} = 0.
\end{eqnarray}
This implies that $p\indd{l}$ and $\rho\indd{l}e\indd{l}$ are constant in space up to the second order perturbation $p\indd{l,(2)}$ and $\left(\rho\indd{l}e\indd{l}\right)\indd{(2)}$.
Furthermore, from the relative velocity equation \eqref{eqn.relvelMS.div.nd} we have the following conditions
\begin{equation}
    \label{eq.condmu}
    \nabla \mu\indd{1,(0)} = \nabla \mu\indd{2,(0)}, \quad \text{and}\quad \nabla \mu\indd{1,(1)} = \nabla \mu\indd{2,(1)}.
\end{equation}
In particular this means that the difference of the chemical potentials $\mu\indd{(0)}$ and $\mu\indd{(1)}$ are constant in space.
Taking these observations into account as well as $\ww_{(0)} = 0$, $\mathcal{O}(M^{-2})$ order terms in the energy equation \eqref{eqn.energyMS.div.nd} reduce to
\begin{equation}
    \label{eq.limit.inte}
    \frac{\pd \rho\indd{(0)} e\indd{(0)}}{\pd t} + \div (\rho\indd{(0)} e\indd{(0)}\vv\indd{(0)} ) +
    p\indd{(0)} \div \vv\indd{(0)} = 0.
\end{equation}
We further assume that the phase pressure relaxation towards a common pressure is faster than the
characteristic time of pressure wave propagation, and one obtains a uniform mixture pressure
$p\indd{(2)} = \vf\indd{1}\rrat{1}p\indd{1,(2)} + \vf\indd{2}\rrat{2}p\indd{2,(2)}= p\indd{2,(2)}$.
This motivates the following assumption on the dynamics of the volume fraction in the low Mach number limit.

\begin{assumption}[Transport of interfaces]\label{Ann.vf}
    In the low Mach number limit we assume
    \begin{equation}\label{eq.Assumption.vf}
        \frac{\pd \vf\indd{1}}{\pd t}  + \vv^{(0)} \cdot \nabla \vf\indd{1} = 0.
    \end{equation}
\end{assumption}
With this assumption, we can rewrite the energy equation at the leading order \eqref{eq.limit.inte} as follows
\begin{equation}
    \label{eq.limit.inte2}
    \vf\indd{1} \frac{\pd (\rho\indd{1}e\indd{1})\indd{(0)}}{\pd t} + \vf\indd{2} \frac{\pd (\rho\indd{2}e\indd{2})\indd{(0)}}{\pd t} + \rho\indd{(0)}e\indd{(0)}\div \vv\indd{(0)} +
    p\indd{(0)} \div \vv\indd{(0)} = 0.
\end{equation}
In fact, \eqref{eq.limit.inte2} can be written as a convex combination with respect to the volume fraction
\begin{equation}
    \label{eq.limit.inte4}
    \vf\indd{1} \left( \frac{\pd (\rho\indd{1,(0)}e\indd{1,(0)})}{\pd t} +
    \left(\rho\indd{1,(0)}e\indd{1,(0)}  + p\indd{1,(0)}\right)\div \vv\indd{(0)}\right)+
    \vf\indd{2} \left(\frac{\pd (\rho\indd{2,(0)}e\indd{2,(0)})}{\pd t} +
    \left(\rho\indd{2,(0)}e\indd{2,(0)} + p\indd{2,(0)}\right)\div \vv\indd{(0)}\right) = 0.
\end{equation}
Since the volume fraction can be arbitrary under the constraint $0 < \vf_l < 1, l = 1,2$, \eqref{eq.limit.inte4} implies
\begin{equation}
    \label{eq.limit.inte5}
    \frac{\pd (\rho\indd{1,(0)}e\indd{1,(0)})}{\pd t} + \left(\rho\indd{1,(0)}e\indd{1,(0)} + p\indd{1,(0)}\right)\div \vv\indd{(0)} = 0 \quad \text{and} \quad \frac{\pd (\rho\indd{2,(0)}e\indd{2,(0)})}{\pd t} + \left(\rho\indd{2,(0)}e\indd{2,(0)} + p\indd{2,(0)}\right)\div \vv\indd{(0)} = 0.
\end{equation}
which is consistent with the limit of single phase flow of the Euler equations, see e.g.
\cite{Dellacherie2010,KlaiMaj1981,Klein1995}.

Then, analogously to \cite{Dellacherie2010} for the case of the Euler equations, we obtain from
\eqref{eq.limit.inte5}, that $\rho\indd{l,(0)}e\indd{l,(0)}$ are constant in space \textit{and} time and consequently we
obtain from \eqref{eq.Mach.expansion.EOS} that also the phase pressures at leading order
$p\indd{l,(0)}$ and $\rho\indd{l,(0)}T\indd{(0)}$ are constant in space and time.
Furthermore, we obtain the divergence free mixture velocity constraint $\div \vv_{(0)} = 0$.
Summarizing, we can formally write the following expansions for the pressure and internal energy
\begin{subequations}
    \begin{align}
        p &= p\indd{(0)} + \mathcal{O}(M^2), \hspace{1.6cm} p\indd{(0)} = \text{constant,} \label{eq.p.exp1}\\
        \rho\indd{l}e\indd{l} & = \rho\indd{l,(0)}e\indd{l,(0)} + \mathcal{O}(M^2),
        \qquad \rho\indd{l,(0)}e\indd{l,(0)}=\text{constant.} \label{eq.rhoe.exp1}
    \end{align}
\end{subequations}
To obtain the expansion for the temperature, we look at the constraint for the chemical potentials.
First, \eqref{eq.rhoe.exp1} implies that $\rho\indd{l,(0)}T\indd{(0)}$ are constant.
Therefore, we can define two constants $\mathcal{E}\indd{l} > 0$ such that
\begin{equation}
    \label{eq.scaling_T}
    T\indd{(0)} \rho\indd{1,(0)} = \mathcal{E}\indd{1}, \quad T\indd{(0)} \rho\indd{2,(0)} = \mathcal{E}\indd{2}.
\end{equation}
Then it follows from $\nabla (\mu\indd{1,(0)} - \mu\indd{2,(0)}) = 0$ that
\begin{equation}
    0 =
    \left(c\indd{v,2}\log\left(\left(\frac{\rho\indd{0,2}}{\mathcal{E}\indd{2}}
    \right)^{\gamma\indd{2}-1}\frac{T\indd{(0)}^{\gamma\indd{2}}}{T\indd{0}}
    \right)-c\indd{v,l}\log\left(\left(\frac{\rho\indd{0,1}}{\mathcal{E}\indd{1}}
    \right)^{\gamma\indd{1}-1}\frac{T\indd{(0)}^{\gamma\indd{1}}}{T\indd{0}}\right)\right)\nabla
    T\indd{(0)}
\end{equation}
and consequently $T\indd{(0)}$ is constant unless both phases coincide.
Since we consider a general case of different phase densities, it follows from \eqref{eq.scaling_T}, that the phase densities $\rho\indd{l,(0)}$ are constant as well.
From relation \eqref{eq.p.exp1} follows $p\indd{(1)} = 0$ thus $\rho\indd{l,(0)}T\indd{(1)}+\rho\indd{l,(1)}T\indd{(0)} = 0$ which implies
\begin{equation}
    \label{eq.cond.T1}
    T\indd{(1)} = - \frac{\rho\indd{l,(1)}}{\rho\indd{l,(0)}}T\indd{(0)}.
\end{equation}
Relation \eqref{eq.condmu} yields $T\indd{(1)}$ is constant and together with \eqref{eq.cond.T1} implies
\begin{equation}
    \label{eq.mu.T.rho.exp1}
    \mu\indd{l} = \mu\indd{l,(0)} + \mathcal{O}(M^2),\quad \rho\indd{l} = \rho\indd{l,(0)} +
    \mathcal{O}(M^2), \quad T = T\indd{(0)} + \mathcal{O}(M^2),\quad
    \{\mu\indd{l,(0)},~\rho\indd{l,(0)}~,T\indd{(0)}\} = \text{constant.}
\end{equation}
Plugging these expansions in \eqref{eqn.contiMS1.div.nd} and \eqref{eqn.contiMS2.div.nd} and single out the first order perturbations for the phase densities, it follows $\vv_{(1)} = 0$ and $\ww_{(1)} = 0$.
Consequently, we have for the friction source term $\tau^{(\ww)} = \mathcal{O}(M^2)$.
We proceed by defining a set of well-prepared initial data.
\begin{definition}[Well-prepared initial data for variable volume fraction]
    \label{def.wp}
    Let $\q \in \mathbb{R}^{4+2d}$ denote the state vector and let both phases be in the same Mach number regime denoted by $M$.
    Let Assumption \ref{Ann.vf} hold.
    Then the set of well-prepared initial data is given as
    \begin{multline}
        \Omega_{wp}^M = \left\lbrace ~\q \in \mathbb{R}^{4+2d}:\quad \rrat{1} p\indu{2}\indd{(k)} = \rrat{1}p\indd{1,(k)}, ~k=0,1,2; \quad \nabla p\indd{1,(0)} = 0, ~\nabla p\indd{2,(0)} = 0;  \quad \nabla p\indd{1,(1)} = 0, ~\nabla p\indd{2,(1)} = 0;\right. \\
        \left. ~\nabla \mu\indd{1,(0)} = \nabla \mu\indd{1,(0)}, ~\nabla \mu\indd{1,(1)} = \nabla \mu\indd{2,(1)}; \quad \div \vv_{(0)} = 0, \quad \vv_{(1)} = 0; \quad \tau^{(\ww)} =\mathcal{O}(M^2) ~\right\rbrace
    \end{multline}
    using the Mach number expansions \eqref{eq.Mach.expansion.phase},\eqref{eq.Mach.expansion.EOS}.
\end{definition}
For well-prepared initial data, we obtain formally for $M \to 0$ and $\tau^{(\ww)} = \mathcal{O}(M^2)$ the following incompressible limit equations with variable volume fraction
\begin{subequations}\label{eq.SHTC.div.limit}
    \begin{eqnarray}
        && \frac{\pd \vf\indd{1}}{\pd t}+\vv\indd{(0)} \cdot \nabla \vf\indd{1} = 0,\quad \rho\indd{(0)} = \vf\indd{1}\rho\indd{1,(0)}+\vf\indd{2} \rho\indd{2,(0)},\label{eqn.alphaMS.div.limit}\\[2mm]
        &&\frac{\pd \vv\indd{(0)}}{\pd t}+
        \vv\indd{(0)}\cdot \nabla \vv\indd{(0)} + \frac{\nabla p\indd{(2)}}{\rho\indd{(0)}} =0, \quad \div \vv_{(0)} = 0,
        \label{eqn.momentumMS.div.limit}
    \end{eqnarray}
\end{subequations}
where $p\indd{(2)}$ is the second order perturbation of the uniform pressure given by $-\div \left(\nabla p\indd{(2)}/\rho\indd{(0)}\right) = \nabla ^2 : \left(\vv_{(0)} \otimes \vv_{(0)}\right)$ acting as the Lagrangian multiplier.
Note that the limit velocity equation is derived applying $\partial_t \vf = - \vv_{(0)} \cdot \nabla \vf$ in the momentum formulation \eqref{eqn.momentumMS.div.nd}.

\section{Numerical scheme}
\label{sec:NumScheme}
Let us write the two-fluid model \eqref{eq.SHTC.div.nd} in the following compact form
\begin{equation}
    \label{eq.short.notation}
    \frac{\pd \q}{\pd t} + \div \bm f(\q)  + \bm B(\q) \cdot \nabla \q = \bm r(\bm q),
\end{equation}
were $\q$ denotes the vector of state variables defined in \eqref{eq.stateVar}, $\bm f$ the flux
function consisting of the conservative terms, $\bm B(\q)$ is the matrix that contains the
non-conservative contributions and $\bm r$ the relaxation source terms acting on the volume
fraction and the relative velocity.

In the following, we construct a numerical scheme for the two fluid single-temperature model
\eqref{eq.short.notation} which is stable independently of the Mach numbers
$M\indu{1}$ and $M\indu{2}$.
This allows to follow the dynamics associated with the flow velocity $\vv$, especially the
transport of the volume fraction that represents the interface between the two phases.
In addition we require the new scheme to be asymptotically preserving (AP), meaning that the
numerical scheme in the singular limit as $M\to 0$ has to be consistent with a discretization of the incompressible
limit equations~\eqref{eq.SHTC.div.limit}.

To achieve this goal, we use an operator splitting approach, dividing the flux $\bm f$ into terms $\bm f\indu{\text{ex}}$ treated
explicitly and $\bm f\indu{\text{im}}$ integrated implicitly.
The components of the non-conservative terms $\bm B$ only involve terms with respect to the velocities $\vv$ and $\ww$ and are treated explicitly.
The resulting implicit system is in general nonlinear due to the nonlinearity of the  flux function
$\bm f\indu{\text{im}}$ and nonlinearity of EOS \eqref{eq.EOS}.
This, however, implies a huge computational overhead since nonlinear solvers would be required to solve large coupled implicit systems.
To reduce computational costs, we construct a linear implicit numerical scheme whose implicit part can be solved by direct or iterative linear solvers.
To avoid that the AP property is lost during the linearization process, we use the so called {\em reference solution (RS) approach} detailed in the subsequent section that has been successfully applied to construct schemes for the Euler equations and isentropic two-phase flows \cite{NoeBisAruLukMun2014,BisLukYel2017,KaiSchSchNoe2017,LukPupTho2022}.

\subsection{Reference solution approach}

In the singular limit as $M\to 0$, the stiffness in the system is mainly connected to the pressure and chemical potential terms in the momentum and relative velocity equation.
Further, these terms are coupled with the evolution equation for the total energy density $\rho E$.
Therefore, to obtain a time step that is dominated by the mixture velocity $\vv$, these terms need to be treated implicitly.
It follows from the EOS \eqref{eq.EOS}, that the mixture pressure $p$ depends linearly on $\rho E$ and we can write
\begin{equation}
    \label{eq.p.expandsion.E}
    p = \left(\phi_p - 1\right)\left(\rho E - \rho E_{\tkin}\right), \quad \text{with} \quad \phi_p(\vf\indd{1},\vf\indd{2},\rho\indd{1},\rho\indd{2}) = \frac{\gamma\indd{1}\vf\indd{1} \rho\indd{1} c\indd{v,1} +\gamma\indd{2} \vf\indd{2}\rho\indd{2} c\indd{v,2}}{\vf\indd{1} \rho\indd{1} c\indd{v,1} + \vf\indd{2}\rho\indd{2} c\indd{v,2}},
\end{equation}
where $E\indd{\tkin}$ contains all kinetic energy contributions
\begin{equation}
    E\indd{\tkin} = \frac{\|\vv\|^2}{2} + \mf\indd{1}\mf\indd{2}\frac{\|\ww\|^2}{2}.
\end{equation}
For the chemical potentials we have a nonlinear dependence on $\rho E$ via the phase entropies $s\indd{l}(\q), l = 1, 2$, given by
\begin{equation}
    \label{eq.mu.expansion.E}
    \mu = \frac{\mu\indd{1}}{M\indd{1}^2} - \frac{\mu\indd{2}}{M\indd{2}^2}= \phi_\mu (\rho E - \rho
    E\indd{\tkin}), \quad \text{where}\quad \phi_\mu(\q) = \frac{\gamma\indd{1} c\indd{v,1} -
        s\indd{1}(\q) - \mathcal{C}^2(\gamma\indd{2} c\indd{v,2} - s\indd{2}(\q))}{\vf\indd{1} \rho\indd{1}
        c\indd{v,1} + \mathcal{C}^2 \vf\indd{2}\rho\indd{2} c\indd{v,2}}
\end{equation}
with $\mathcal{C}$ being the ratio of Mach numbers defined in \eqref{eq.Mach_C}.

Note that for single phase flow, $\phi_p = \gamma -1$ is constant and $\phi_\mu = 0$.
Consequently, formulations \eqref{eq.p.expandsion.E}, \eqref{eq.mu.expansion.E} reduce to the Euler case, studied in \cite{BosRusSca2018}, and is consistent with single phase flow.
Therefore, we linearize only the difference of the chemical potentials $\mu$ with respect to $\rho E$ around a reference state $\q\indd{\RS}$ as follows
\begin{equation}
    \mu = \mu\indd{\RS} + \left(\frac{\pd \mu}{\pd (\rho E)}\right)\indd{\RS} \left(\rho E - (\rho E)\indd{\RS}\right) + \mathcal{O}\left(\left(\rho E - (\rho E)\indd{\RS}\right)^2\right),
\end{equation}
where
\begin{equation}
    \label{eq.dmu.exp}
    \frac{\pd \mu}{\pd (\rho E)} = \frac{c\indd{v,1} (\gamma\indd{1} - 1) - s\indd{1}(\rho\indd{1},T) -\mathcal{C}^2(c\indd{v,2}(\gamma\indd{2} -1)- s\indd{2}(\rho\indd{2},T))}{\vf\indd{1} \rho\indd{1} c\indd{v,1} + \mathcal{C}^2\vf\indd{2}\rho\indd{2} c\indd{v,2}} = \mathcal{O}(1)
\end{equation}
for all $\mathcal{C} = \frac{M\indd{1}}{M\indu{2}}$.
The reference state is set to
\begin{equation}\label{def.RSstates}
    \q\indd{\RS} = (\vf\indd{1}, \vf\indd{1}\rho\indd{1,(0)}, \vf\indd{2}\rho\indd{2,(0)}, \rho\indd{(0)}\vv, \ww, (\rho E)\indd{\RS})^T, \quad (\rho E)\indd{\RS} = \rho\indd{(0)}e\indd{(0)} + \rho\indd{(0)}E\indd{\tkin},
\end{equation}
where $\rho\indd{l,(0)}, \rho\indd{l,(0)}e\indd{l,(0)}$ are constant leading order states from \eqref{eq.mu.T.rho.exp1} and $$\rho\indd{(0)} = \vf\indd{1} \rho\indd{1,(0)} + \vf\indd{2} \rho \indd{2,(0)}, \quad \rho\indd{(0)}e\indd{(0)} = \vf\indd{1} \frac{\rho\indd{1,(0)}e\indd{1,(0)}}{M\indd{1}^2} + \vf\indd{2}\frac{\rho\indd{2,(0)} e\indd{2,(0)}}{M\indd{2}^2}.$$
We split $\mu$ into a part that is linear in $\rho E$ given by
\begin{equation}
    \hat \mu = \hat\mu\indd{\RS} + \left(\frac{\pd \mu}{\pd (\rho E)}\right)\indd{\RS} \rho E, \quad \hat \mu\indd{\RS} =  \mu(\q\indd{\RS}) - \left(\frac{\pd \mu}{\pd (\rho E)}\right)\indd{\RS} \times (\rho E)\indd{\RS}
\end{equation}
and a nonlinear part
\begin{equation}
    \bar \mu = \mu - \hat \mu = \mathcal{O}\left(\left(\rho E - (\rho E)\indd{\RS}\right)^2\right).
\end{equation}
Note that if considering well-prepared initial data, $\bar \mu$ is of order $M\indd{l}^2, l=1,2$.
This can be seen by multiplying $\mu = \frac{\mu\indd{1}}{M\indd{1}^2} - \frac{\mu\indd{2}}{M\indd{2}^2}$, without loss of generality, by $M\indd{1}^2$.
Then we obtain
\begin{equation}
    \mu\indd{1} - \mathcal{C}^2 \mu\indd{2} = \mu\indd{1,\RS} - \mathcal{C}^2 \mu\indd{2,\RS} + \left( \frac{\pd \mu}{\pd (\rho E)}\right)\indd{\RS} \times \left(\rho \check{E} - (\rho \check{E})\indd{\RS}\right) + \mathcal{O}\left((\rho \check{E} - (\rho \check {E})\indd{\RS})^2\right),
\end{equation}
where
\begin{equation}
    \label{eq.scaled_total_energy}
    \rho\check{E} = \vf\indd{1}\rho\indd{1}e\indd{1} + \mathcal{C}^2\vf\indd{2}\rho\indd{2}e\indd{2} + M\indd{1}^2 \rho E\indd{\tkin}.
\end{equation}
For the difference between the total energy density and its reference value we obtain
\begin{equation}
    \rho\check{E} - (\rho\check{E})\indd{\RS} = M\indd{1}^2\vf\indd{1}\left(\rho\indd{1}e\indd{1}\right)\indd{(2)} + (M\indd{2}\mathcal{C})^2\vf\indd{2}\left(\rho\indd{2}e\indd{2}\right)\indd{(2)} + M\indd{1}^2\rho\indd{(2)}E\indd{\tkin}  = M\indd{1}^2\left(\mf\indd{1}e\indd{1,(2)} + \mf\indd{2}e\indd{2,(2)} +\rho\indd{(2)}E\indd{\tkin}\right)
\end{equation}
with $\rho\indd{(2)} = M\indd{1}^2\vf\indd{1}\rho\indd{1,(2)} + M\indd{2}^2\vf\indd{2}\rho\indd{2,(2)}$.
Therefore, it follows
\begin{equation}
    \bar \mu = \frac{\mu\indd{1} - \mathcal{C}^2 \mu\indd{2}}{M\indd{1}^2} - \hat \mu = \mathcal{O}\left(\frac{(\rho \check{E} - \rho \check {E}\indd{\RS})^2}{M\indd{1}^2}\right) = \mathcal{O}(M\indd{1}^2).
\end{equation}
Consequently, the nonlinear term $\bar \mu$ vanishes as $M\indd{1}$ tends to 0 and can be treated explicitly without imposing a sever time step restriction.
Note that for compressible flow $M\indd{1} \approx 1$, these terms are important to obtain the correct wave speeds and cannot be neglected.

Taking these considerations into account, the following subsystem will be treated explicitly
\begin{subequations}\label{eq.SHTC.div.nd.expl}
    \begin{eqnarray}
        && \frac{\pd \vf\indd{1}}{\pd t}+\vv \cdot \nabla \vf\indd{1} =0,\label{eqn.alphaMS.div.nd.expl}\\[2mm]
        && \frac{\pd (\vf\indd{1}\rho\indd{1})}{\pd t} +\div (\vf\indd{1}\rho\indd{1} \vv\indd{1}) = 0,\label{eqn.contiMS1.div.nd.expl}\\[2mm]
        && \frac{\pd (\vf\indd{2}\rho\indd{2})}{\pd t}+\div (\vf\indd{2}\rho\indd{2} \vv\indd{2}) = 0,\label{eqn.contiMS2.div.nd.expl}\\[2mm]
        &&\frac{\pd (\rho \vv)}{\pd t}+\div
        \left(\rho \vv \otimes \vv\right)=0,
        \label{eqn.momentumMS.div.nd.expl}\\[2mm]
        &&\frac{\pd \ww}{\pd t}+\div\left(\left[ \ww \cdot \vv + \left(1-2\mf\indd{1}\right) \frac{\|\ww\|^2}{2}\right]\Id\right) + ( \nabla \times \ww ) \times \vv=0,\label{eqn.relvelMS.div.nd.expl}\\[2mm]
        &&
        \frac{\pd (\rho E)}{\pd t}+ \div \left(\rho ~\left[ \ww \cdot \vv+ \left(1-2\mf\indd{1}\right) \frac{\|\ww\|^2}{2}\right] \mf\indd{1}\mf\indd{2}\ww \right)=0.\label{eqn.energyMS.div.nd.expl}
    \end{eqnarray}
\end{subequations}
Written in compact notation, we have
\begin{equation}
    \label{eq.subsys.ex}
    \frac{\pd \q}{\pd t} + \div \bm f^{\text{ex}}(\q) + \bm B(\q) \nabla \q = 0.
\end{equation}
Subsystem \eqref{eq.SHTC.div.nd.expl} is weakly hyperbolic, since it lacks one linearly independent eigenvector for the characteristic speed $\lambda_1^\vv$.
The complete list of characteristic speeds is given by
\begin{equation}
    \label{eq.CharSpeeds}
    \lambda^0 = 0, \quad
    \lambda^\vv_1 = \vv \cdot \nn ~(8 \times), \quad
    \lambda^\vv_2 = (\vv + (1-2 \mf\indd{1}) \ww) \cdot \nn.
\end{equation}
Applying, e.g. the Rusanov numerical fluxes, a numerical solution of the weakly hyperbolic system \eqref{eq.SHTC.div.nd.expl} can be obtained.
Further, rewriting pressures and chemical potentials in terms of $\rho E$ and using the decomposition $$\mu = \hat \mu + \bar \mu \quad \text{and }\quad\partial \mu\indd{RS} = \left(\frac{\pd \mu}{\pd(\rho E)}\right)\indd{\RS},$$ the implicitly treated subsystem is given by
\begin{subequations}\label{eq.SHTC.div.nd.impl}
    \begin{eqnarray}
        && \frac{\pd \vf\indd{1}}{\pd t} =-\frac{1}{\tau^{(\vf)}\rho}\left(\rrat{1}\frac{p\indd{1}}{M\indd{1}^2}-\rrat{2}\frac{p\indd{2}}{M\indd{2}^2}\right),\label{eqn.alphaMS.div.nd.impl}\\[2mm]
        && \frac{\pd (\vf\indd{1}\rho\indd{1})}{\pd t} = 0,\label{eqn.contiMS1.div.nd.impl}\\[2mm]
        && \frac{\pd (\vf\indd{2}\rho\indd{2})}{\pd t} = 0,\label{eqn.contiMS2.div.nd.impl}\\[2mm]
        &&\frac{\pd (\rho \vv)}{\pd t}+\div
        \left((\phi_p - 1)\left( \rho E - \rho E\indd{\tkin}\right)\Id + \rho \mf\indd{1}\mf\indd{2} \ww \otimes \ww \right)=0, \label{eqn.momentumMS.div.nd.impl}\\[2mm]
        &&\frac{\pd \ww}{\pd t}+\div\left(\left[ \partial\mu\indd{\RS}~\rho E + \hat\mu\indd{\RS} + \bar \mu \right]\Id\right)=-\frac{
            \mf\indd{1}\mf\indd{2}\ww }{\tau^{(w)}},\label{eqn.relvelMS.div.nd.impl}\\[2mm]
        &&
        \frac{\pd (\rho E)}{\pd t}+
        \div \left(\rho \vv \left(\frac{\rho E + p}{\rho}\right)  + \mu~ \rho \mf\indd{1}\mf\indd{2}\ww \right)=0\label{eqn.energyMS.div.nd.impl}
    \end{eqnarray}
\end{subequations}
which yields the corresponding compact form
\begin{equation}
    \label{eq.subsys.im}
    \frac{\pd \q}{\pd t} + \div \bm f^{\text{im}}(\q) = \bm r(\q).
\end{equation}
In the following, we construct an IMEX scheme for two subsystems \eqref{eq.subsys.ex} and \eqref{eq.subsys.im}.
We start with the time semi-discrete scheme.

\subsection{Time semi-discrete scheme}\label{sec.SemiDiscrScheme}
Let the time interval $(0,T_f)$ be discretized by $t^n = n \Delta t$, where $\Delta t$ denotes the time step subject to a time step restriction based on a CFL condition given by
\begin{equation}
    \label{eqn.CFL.cond}
    \Delta t \leq \nu \frac{\Delta x}{\underset{{\bm x \in \Omega} }{\max }~(|\lambda_1^\vv(\bm x, t^n)|, |\lambda_2^\vv(\bm x, t^n)|)}.
\end{equation}
Therein, $\lambda_1^\vv, ~\lambda_2^\vv$ are the characteristic speeds of the explicit subsystem \eqref{eq.CharSpeeds} evaluated at time level $t^n$.
For a first order scheme in time, we apply the forward Euler method for the explicit subsystem \eqref{eq.subsys.ex}
\begin{equation}
    \label{eq.time-semi.expl}
    \q^{(1)} = \q^n - \Delta t \div \bm f^{ex}(\q^n) - \Delta t \bm B(\q^n) \cdot \nabla \q^n
\end{equation}
and a backward Euler method for the implicit subsystem \eqref{eq.SHTC.div.nd.impl}.
We find that there are still some nonlinear terms present yielding a nonlinear coupled system.
Extending the approach from \cite{BosRusSca2018} for the Euler equations, we linearize certain flux terms in time yielding the following time discretization
\begin{subequations}
    \label{eq.implicit_sub_part}
    \begin{eqnarray}
        &&(\vf\indd{1}\rho\indd{1})^{n+1} = (\vf\indd{1}\rho\indd{1})^{(1)},\\[2mm]
        &&(\vf\indd{2}\rho\indd{2})^{n+1} = (\vf\indd{2}\rho\indd{2})^{(1)},\\[2mm]
        &&(\rho \vv)^{n+1} = \rho\vv^{(1)} - \Delta t\nabla
        \left((\phi_p^n - 1) \rho E^{n+1} + \hat p^{n}\right) - \Delta t \div \left(\left(\rho \mf\indd{1}\mf\indd{2}\right)^{n+1}\ww^{n} \otimes \ww^{n} \right),\\[2mm]
        &&\ww^{n+1}=\ww^{(1)}-\Delta t\nabla\left(\pd\mu^n_\RS (\rho E)^{n+1} + \hat\mu\indd{\RS}^{(1)} + \bar \mu^n \right)- \frac{ \Delta t
        }{\tau^{(w)}}\left(\mf\indd{1}\mf\indd{2}\right)^{n+1}\ww^{n+1},\\[2mm]
        &&
        (\rho E)^{n+1} = (\rho E)^{(1)}- \Delta t
        \div \left((\rho\vv)^{n+1} \frac{\left(\rho E + p\right)^n}{\rho^{n+1}} + \mu^n (\rho\mf\indd{1}\mf\indd{2})^{n+1}\ww^{n+1} \right) \label{eq.implicit_sub_part_energy},
    \end{eqnarray}
\end{subequations}
where $\hat{p}^{n} = -(\phi_p^n - 1) \rho E_\tkin^n$.
Rewriting the relative velocity equations with $\rho^{n+1} = (\vf\indd{1}\rho\indd{1})^{n+1} + (\vf\indd{2}\rho\indd{2})^{n+1}$ implies
\begin{equation}
    \label{eq.time-sem-rel-vel}
    \ww^{n+1} = \left(\frac{\tau^{(w)}}{\tau^{(w)} + \Delta t
        \left(\mf\indd{1}\mf\indd{2}\right)^{n+1} }\right)\ww^{(1)}-\left(\frac{\Delta t~\tau^{(w)}}{\tau^{(w)} + \Delta t
        \left(\mf\indd{1}\mf\indd{2}\right)^{n+1} }\right)\nabla\left(\pd\mu^n_\RS (\rho E)^{n+1} + \hat\mu\indd{\RS}^{n} + \bar \mu^n \right).
\end{equation}
Substituting the relative velocity and momentum in the total energy equation yields a linear implicit equation for the total energy given by
\begin{align}
    \label{eq.impl.Etot}
    \begin{split}
        &(\rho E)^{n+1}  - \Delta t^2
        \div \left(\frac{\left(\rho E + p\right)^n}{\rho^{n+1}} \nabla
        \left((\phi_p^n - 1) (\rho E)^{n+1}\right)\right) - \Delta t^2 \div\left( \left(\frac{\tau^{(w)} \mu^n(\rho\mf\indd{1}\mf\indd{2})^{n+1}}{\tau^{(w)} + \Delta t
            \left(\mf\indd{1}\mf\indd{2}\right)^{n+1} }\right)\nabla\left(\pd\mu^n_\RS (\rho E)^{n+1} \right) \right) \\
        = ~&(\rho E)^{(1)} -\Delta t
        \div \left(\frac{\left(\rho E + p\right)^n}{\rho^{n+1}}\left(\rho\vv^{(1)} + \Delta t \nabla \left((\phi_p-1)\rho E\indd{\tkin}\right)^n -\Delta t\div \left(\left(\rho \mf\indd{1}\mf\indd{2}\right)^{n+1}\ww^{n} \otimes \ww^{n} \right)\right)\right) \\
        & -  \Delta t \div\left(\left(\frac{\tau^{(w)} \mu^n(\rho\mf\indd{1}\mf\indd{2})^{n+1}}{\tau^{(w)} + \Delta t
            \left(\mf\indd{1}\mf\indd{2}\right)^{n+1} }\right)\ww^{(1)}-\left(\frac{\Delta t~\tau^{(w)} \mu^n(\rho\mf\indd{1}\mf\indd{2})^{n+1}}{\tau^{(w)} + \Delta t
            \left(\mf\indd{1}\mf\indd{2}\right)^{n+1} }\right)\nabla\left(\hat\mu\indd{\RS}^{n} + \bar \mu^n \right)\right).
    \end{split}
\end{align}
Having obtained the total energy $(\rho E)^{n+1}$ we can successively update the relative velocity
\begin{equation}
    \label{eq.impl.w}
    \ww^{n+1} = \left(\frac{\tau^{(w)}}{\tau^{(w)} + \Delta t
        \left(\mf\indd{1}\mf\indd{2}\right)^{n+1} }\right)\left(\ww^{(1)}-\Delta t \nabla \mu^{n+1} \right)
\end{equation}
and the momentum
\begin{equation}
    \label{eq.impl.mom}
    (\rho \vv)^{n+1} =  \rho\vv^{(1)} - \Delta t \nabla
    p^{n+1} - \Delta t \div \left(\left(\rho \mf\indd{1}\mf\indd{2}\right)^{n+1}\ww^{n+1} \otimes \ww^{n+1} \right).
\end{equation}
Finally, the volume fraction at the next time level is obtained from the pressure relaxation
\begin{equation}
    \frac{\pd \vf\indd{1}}{\pd t} =-\frac{1}{\tau^{(\vf)}\rho}\left(\frac{p\indd{1}}{M\indd{1}^2}-\frac{p\indd{2}}{M\indd{2}^2}\right).
\end{equation}
Rewriting the source term in terms of the state variables $\q$, we find
\begin{eqnarray}
    \frac{\pd \vf\indd{1}}{\pd t} = - \frac{1}{\tau^{(\vf)}}\left(\frac{(\gamma\indd{1}-1)c\indd{v,1}\mf\indd{1}}{\vf\indd{1}} - \frac{(\gamma\indd{2}-1)c\indd{v,2}\mf\indd{2}}{1-\vf\indd{1}}\right)\frac{\rho E - \rho E\indd{\tkin}}{\vf\indd{1}\rho\indd{1}c\indd{v,1} + \vf\indd{2}\rho\indd{2}c\indd{v,2}}
\end{eqnarray}
which is a nonlinear ordinary differential equation in $\vf\indd{1}$ and can be solved implicitly applying the backward Euler scheme and the Newton algorithm to solve the nonlinear implicit system which concludes the time semi-discrete scheme.
We proceed with the construction of the fully discrete scheme in the next section.

\subsection{Fully discrete scheme}\label{sec.FullDiscrScheme}
In time, we set as before $t^{n+1} = t^n + \Delta t$, where $\Delta t$ obeys the CFL condition \eqref{eqn.CFL.cond}.
In space, we consider a two-dimensional computational domain $\Omega$ divided into cells $C_I = [x_{1,i-1/2}, x_{1,i+1/2}]\times[x_{2,j-1/2}, x_{2,j+1/2}]$ with $\bm x = (x_1, x_2)^T$.
The common edge between two neighboring cells $\Omega_I$ and $\Omega_J$ is denoted by $\partial\Omega_{IJ}$ and the set of neighbors of $\Omega_I$ associated with the unit normal vector pointing from the cell $\Omega_I$ to $\Omega_J$ given by $\nn_{IJ}$ is denoted by $\mathcal{N}_I$.
We consider a uniform mesh size $\Delta x_1, \Delta x_2$ in each direction and the barycenter of $C_I$ is denoted by $\bm x_I = (i\Delta x_1, j\Delta x_2)$ for $i,j = 1, \dots, N$.
We use a finite volume framework, where the solution on the cell $C_I$ at time $t^n$ is approximated by the average given by
\begin{equation}
    \q_I^n \approx \frac{1}{|\Omega_I|}\int_{\Omega_I} \q(\bm x, t^n) ~d \bm x.
\end{equation}
A fully discrete finite volume (FV) method for \eqref{eq.subsys.ex} reads
\begin{equation}
    \label{eq.full.update.ex}
    \q_I^{(1)} = \q_I^{n} - \Delta t \sum_{K \in \mathcal{N}_I} \frac{|\partial \Omega_{IK}|}{|\Omega_I|} \left(\bm F^{\text{ex}}(\q_I^n,\q_K^n) \cdot \nn_{IK} + \bm D(\bm q_I^n,\bm q_K^n) \cdot \nn_{IK} \right),
\end{equation}
using a Rusanov numerical flux
\begin{equation}
    \label{eq.Rusanov}
    \bm F^{\text{ex}}(\q_I^n,\q_K^n) \cdot \nn_{IK} = \frac{1}{2}\left(\bm f^{\text{ex}}(\q_I^n) +
    \bm f^{\text{ex}}(\q_K^n) \right) \cdot \nn_{IK} - s_{IK} ~\bm{I}~ (\q_K^n - \q_I^n),
\end{equation}
where $s = \max_k(|\lambda_k(\bm q_I)|,|\lambda(\bm q_I)|)$ denotes the maximum eigenvalue at the interface $\partial\Omega_{IK}$.
The non-conservative product is approximated in the following way
\begin{equation}\label{eq.NonConsApprox}
    \bm D(\bm q_I^n,\bm q_K^n) \cdot \nn_{IK} = \frac{1}{2} \bm B\left(\tilde \q^n\right) \cdot (\q_K^n - \q_I^n), \quad \tilde\q = \frac{1}{2}\left(\q_K^n + \q_I^n\right).
\end{equation}
The implicit elliptic equation for the total energy \eqref{eq.impl.Etot} is based on centered finite difference approximation for the space discretization and can be formulated on cell $C_I$ in the following way
\begin{align}
    \label{eq.implicit_energy_fully_disc}
    \begin{split}
    \left(\rho E\right)_I^{n+1} - \Delta t^2 \left(\mathcal{L}_I^n \left(\rho E\right)_I^{n+1} + \mathcal{K}_I^n \left(\rho E\right)_I^{n+1} \right) =& \left(\rho E\right)_I^{(1)} - \Delta t \sum_{K \in \mathcal{N}_I} \frac{|\delta \Omega_{IK}|}{|\Omega_I|} \mathcal{F}(\q^{(1)}_I,\q_K^{(1)}) \cdot \nn_{IK} \\
    &- \Delta t^2 \mathcal{L}_I^n \left(\rho E_\tkin^n \right) - \Delta t^2 \mathcal{K}_I^n \left((\rho E)_{\RS}^{n} - (\pd \mu_\RS^n)^{-1}\bar \mu^n\right).
    \end{split}
\end{align}
Here, the weighted Laplacians are discretized as follows
\begin{equation}
    \label{eq.Elliptic_fully_disc}
    \mathcal{L}_I \left(\rho E\right)_I = \sum_{K \in \mathcal{N}_I} \frac{|\partial \Omega_{IK}|}{|\Omega_I|} G_1(\q_I,\q_K) \left[H_1 (\rho E)\right](\q_I,\q_K), \quad \mathcal{K}_I \left(\rho E\right)_I = \sum_{K \in \mathcal{N}_I} \frac{|\partial \Omega_{IK}|}{|\Omega_I|} G_2(\q_I,\q_K) \left[H_2(\rho E)\right](\q_I,\q_K)
\end{equation}
with
\begin{equation}
    G_k(\q_I,\q_K) = \frac{1}{2}\left(g_k(\q_I) + g_k(\q_K)\right), \quad H_k(\q_I,\q_K) = \frac{|\partial \Omega_{IK}|}{|\Omega_I|}\left(h_k(\q_I) - h_k(\q_K)\right), \quad k = 1, 2
\end{equation}
where
\begin{equation}
    g_1 = \frac{(\rho E + p)^n}{\rho^{n+1}}, \quad h_1 = \phi_p^n - 1, \quad g_2 = \frac{\tau^{(w)} \mu^n(\rho\mf\indd{1}\mf\indd{2})^{n+1}}{\tau^{(w)} + \Delta t
        \left(\mf\indd{1}\mf\indd{2}\right)^{n+1} }, \quad h_2 = \pd\mu^n_\RS.
\end{equation}
The divergence terms are approximated as
\begin{equation}
    \label{eq.F.centred}
    \mathcal{F}(\q_I^{(1)},\q_K^{(1)}) \cdot \nn_{IK} =  \frac{1}{2}\left( g_1(\q_I) \rho \vv^{(1)\ast}_I + g_1(\q_K) \rho \vv^{(1)\ast}_K \right) \cdot \nn_{IK},
\end{equation}
where $\rho \vv^{(1)\ast}$ contains the relative velocity flux component $\left((\rho\mf_1\mf_2)^{n+1} \ww^n \otimes \ww^n\right)$ with centered differences analogously to \eqref{eq.F.centred}.
The coefficient matrix resulting from the linear equation \eqref{eq.Elliptic_fully_disc} is strictly diagonal dominant.
Therefore, the linear system of equations has a unique solution independent of the Mach number regime.
Numerically it is solved by a preconditioned linear iterative solver GMRES provided by PetSc
\cite{petsc-user-ref}.

Once the energy $\rho E$ is computed at $t^{n+1}$, the relative velocity \eqref{eq.impl.w} and momentum \eqref{eq.impl.mom} are updated consecutively by the FV method
\begin{subequations}
\begin{align}
        \ww_I^{n+1} &= \ww_I^{(1)} - \Delta t \sum_{K \in \mathcal{N}_I} \frac{|\partial \Omega_{IK}|}{|\Omega_I|} \bm F^{\text{im}}_{(\ww)}(\q_I^{n+1},\q_K^{n+1}) \cdot \nn_{IK} + \Delta t\bm r_{(\ww)}(\q_I^{n+1}),\label{eq.impl.upd.w}\\
        \rho\vv_I^{n+1} &= \rho\vv_I^{(1)} - \Delta t \sum_{K \in \mathcal{N}_I} \frac{|\partial \Omega_{IK}|}{|\Omega_I|} \bm F^{\text{im}}_{(\rho\vv)}(\q_I^{n+1},\q_K^{n+1}) \cdot \nn_{IK}. \label{eq.impl.upd.rhou}
\end{align}
\end{subequations}
The numerical flux $\bm F^{\text{im}}$ is constructed by the finite difference approximation defined analogously as in \eqref{eq.F.centred} based on the implicit flux $\bm f^{\text{im}}$.
The update of the volume fraction is approximated by the backward Euler method
\begin{equation}
    \label{eq.impl.vf.relax}
    \vf_I^{n+1} = \vf_I^{(1)} + \Delta t \bm r_{(\vf)}(\q_I^{n+1}),
\end{equation}
where $\bm r_{(\vf)}$ denotes the pressure relaxation source term.
Update \eqref{eq.impl.vf.relax} approximates on each cell the corresponding  ordinary differential equation (ODE) independently.
To solve the nonlinear implicit system arising from the backward Euler discretization of the ODE, a Newton algorithm is applied.

The steps of the reference solution implicit-explicit finite volume (RS-IMEX FV) scheme can be summarized as follows:
    \begin{enumerate}
        \item Compute the explicit update $\q_I^{(1)}$ given by \eqref{eq.full.update.ex} based on the advective terms $\bm f^\text{ex}$ and $\bm B$ with the numerical flux \eqref{eq.Rusanov} and the approximation of the non-conservative terms \eqref{eq.NonConsApprox} under the material CFL condition \eqref{eqn.CFL.cond}.
        \item Compute the implicit update $\q_I^{n+1}$ given by the following consecutive steps:
        \begin{enumerate}
            \item Solve the linear implicit equation \eqref{eq.implicit_energy_fully_disc} with centred elliptic operators \eqref{eq.Elliptic_fully_disc} for the total energy $(\rho E)_I^{n+1}$ based on a linearization using reference states $\q_\RS^n$ defined in \eqref{def.RSstates}.
            \item Compute first the update of relative velocity $\ww_I^{n+1}$ given in
            \eqref{eq.impl.upd.w} and then the momentum $\rho\vv_I^{n+1}$ given in
            \eqref{eq.impl.upd.rhou} using the respective full nonlinear flux components given in
            $\bm f^\text{im}$ discretized with centred numerical fluxes. Due to the knowledge of
            $(\rho E)_I^{n+1}$ and the consecutive execution, both updates can be done explicitly.
            \item Solve on each cell the nonlinear implicit system \eqref{eq.impl.vf.relax}  for the volume fraction $\vf_I^{n+1}$ arising from the implicit treatment of the pressure relaxation process $\bm r_{(\vf)}$ using a Newton algorithm.
        \end{enumerate}
    \end{enumerate}

\subsection{ Higher order extension}

The above procedure fits in the framework of IMEX Runge Kutta methods, using a forward Euler scheme for the explicit and an backward Euler scheme for the implicit subsystems.
The corresponding Butcher tableaux are given in Table \ref{tab.Butcher_IMEX12}.
\begin{table}
    \renewcommand{\arraystretch}{1.25}
    \centering
    \begin{tabular}{c || c c c}
        IMEX-RK Scheme & \begin{tabular}{c|c}
            $d$ & $A$ \\
            \hline
            & $b^T$ \\
        \end{tabular} & \begin{tabular}{c|c}
            $\tilde d$ & $\tilde A$ \\
            \hline
            & $\tilde b^T$ \\
        \end{tabular}
        &
        \\[5mm]
        \hline \hline \\[-5mm]
        Backward/Forward Euler & \begin{tabular}{c|c}
            1 & 1 \\
            \hline
            & 1 \\
        \end{tabular} & \begin{tabular}{c|c}
            0 & 0 \\
            \hline
            & 1 \\
        \end{tabular} &\\[5mm]
        \hline\hline \\[-5mm]
        ARS(2,2,2) & \begin{tabular}{c|ccc}
            0 & 0 & 0 & 0 \\
            $\gamma$ & 0 &$\gamma$ & 0 \\
            1 & 0& $1-\gamma$ & $\gamma$ \\
            \hline
            & 0& $1-\gamma$ & $\gamma$ \\
        \end{tabular} & {\color{blue}\begin{tabular}{c|ccc}
            0 & 0 & 0 & 0 \\
             $\gamma$& $\gamma$ & 0 & 0 \\
            1 & $\delta$ & $1 - \delta$ & 0 \\
            \hline
            & $\delta$ & $1 - \delta$ & 0\\
        \end{tabular}}
        & $\gamma = 1 - \frac{1}{\sqrt{2}}$, $\delta = 1 - \frac{1}{2\gamma}$
    \end{tabular}
    \caption{Butcher tableaux of the first and second order scheme.}
    \label{tab.Butcher_IMEX12}
\end{table}
For an $s$ stage IMEX Runge Kutta method, the Butcher tableaux $(A,b,d)$ for the implicit and $(\tilde A, \tilde b, \tilde d)$ for the explicit parts are given by
\begin{equation}
    A =
    \begin{bmatrix}
        a_{11} & \cdots & 0 \\
        \vdots & \ddots & \vdots \\
        a_{s1} & \cdots & a_{ss}
    \end{bmatrix},
    \quad
    b = \begin{bmatrix}
        b_1 \\ \vdots \\ b_s
    \end{bmatrix}
    \quad
    d = \begin{bmatrix}
        d_1 \\ \vdots \\ d_s
    \end{bmatrix},
    \quad
    \tilde A =
    \begin{bmatrix}
        \tilde a_{11} & \cdots & 0 \\
        \vdots & \ddots & \vdots \\
        \tilde a_{s1} & \cdots & \tilde a_{ss}
    \end{bmatrix},
    \quad
    \tilde b = \begin{bmatrix}
        \tilde  b_1 \\ \vdots \\ \tilde b_s
    \end{bmatrix}
    \quad
    \tilde d = \begin{bmatrix}
        \tilde d_1 \\ \vdots \\ \tilde d_s
    \end{bmatrix}.
\end{equation}
Consequently, the IMEX method can be written as
\begin{equation}
    \label{eq.IMEX.final}
    \q_I^{n+1} = \q_I^{n} - \Delta t \sum_{k=1}^s {\color{blue}\tilde b_k}\left(\div \bm f^{\text{ex}}(\q^{(k)}) + \bm B(\q^{(k)})\cdot \nabla \q^{(k)}\right) + {\color{blue}b_k}\left(\div \bm f^{\text{im}}(\q^{(k)}) - \bm r(\q^{(k)})\right),
\end{equation}
with the stages {\color{blue} evaluated at time $t^{(k)} = t^n + d_k \Delta t$}
\begin{equation}
    \label{eq.IMEX.stage}
    \q_I^{(k)} = \q_I^{n} - \Delta t \sum_{i=1}^{k-1} \tilde a_{ki}\left(\div \bm f^{\text{ex}}(\q^{(i)}) + \bm B(\q^{(i)})\cdot \nabla \q^{(k)}\right) - \Delta t \sum_{i=1}^{k} a_{ki}\left(\div \bm f^{\text{im}}(\q^{(i)}) - \bm r(\q^{(i)})\right).
\end{equation}
To be consistent with the asymptotic limit, we apply globally stiffly accurate (GSA) IMEX-RK methods like for the time discretization.
For the first order method forward/backward Euler method is applied, for the second order method ARS(2,2,2) is used, see Table \ref{tab.Butcher_IMEX12}.
A second order method in space is achieved by a second order reconstruction with a minmod limiter in the Rusanov numerical flux \eqref{eq.Rusanov}.
In the implicit part, central finite differences yield second order accuracy.
Note that for GSA Runge Kutta methods, the final update \eqref{eq.IMEX.final} coincides with the last computational stage \eqref{eq.IMEX.stage} and thus does not need to be performed.

\section{Asymptotic preserving property}\label{sec:APanalysis}

Motivated by the analysis in Section \ref{sec.WellPrepLowMachLimits}, we consider the case $M\indd{1} = M\indd{2} = M \ll 1$.
For the cases $M\indd{1}\neq M\indd{2}$ and $M\indd{1} \approx 1 \geq M\indd{2} > 0$ we refer to the study of the isentropic case performed in \cite{LukPupTho2022}.
The principle idea is the same and the proof can be performed along the lines presented in \cite{LukPupTho2022} combined with the analysis for the case $M\indd{1} = M\indd{2}$ presented here.

Applying an analogous asymptotic analysis as in Section \ref{sec.WellPrepLowMachLimits} on the semi-discrete scheme consisting of \eqref{eq.time-semi.expl}, \eqref{eq.impl.Etot}, \eqref{eq.impl.w}, \eqref{eq.impl.mom}, \eqref{eq.impl.vf.relax}  and using well-prepared initial data as defined in Definition~\ref{def.wp}, we can prove the asymptotic preserving (AP) property for the semi-discrete scheme.

\begin{theorem}[Asymptotic preserving property]
    \label{theo:AP}
    The first order RS-IMEX FV scheme consisting of the explicit part \eqref{eq.time-semi.expl}, the linear implicit elliptic system \eqref{eq.impl.Etot} and the implicit updates \eqref{eq.impl.w}, \eqref{eq.impl.mom} and \eqref{eq.impl.vf.relax} is asymptotic preserving up to $\mathcal{O}(\Delta t)$.
    More precisely, for well-prepared initial data $\q^0 \in \Omega_{wp}^M$ the RS-IMEX FV scheme yields a consistent approximation of limit equations \eqref{eq.SHTC.div.limit} up to $\mathcal{O}(\Delta t)$.
\end{theorem}

\begin{proof}
    For the proof of the AP property we refer a reader to Appendix \ref{App:APproof}.
\end{proof}

We want to point out, that the $\mathcal{O}(\Delta t)$ errors arising in the velocity equation and divergence free constraint at leading order are due to the non-constant volume fraction in the low Mach number limit.
For single phase flow, i.e. the Euler equations, or homogeneous mixtures, i.e. constant $\vf$, we obtain a stronger result of $\div \vv^{n+1} = \mathcal{O}(M^2)$.

\begin{corollary}[Asymptotic preserving property for constant volume fraction]
    For constant volume fraction or single phase flows,
    we have
    \begin{equation}
        \label{eq.corol.1}
        \vv_{(0)}^{n+1} = \vv_{(0)}^{(1)} - \Delta t \vv_{(0)}^n \cdot \nabla \vv_{(0)}^n - \Delta t \frac{\nabla p_{(2)}^{n+1}}{\rho^{n+1}}.
    \end{equation}
    Moreover, the energy update
    \begin{equation}
        (\rho e)^{n+1}_{(0)} = (\rho e)^n_{(0)} - \Delta t ((\rho e + p )_{(0)}^n \vv_{(0)}^{n+1})
    \end{equation}
    and $T_{(0)}^{n+1} = T_\RS + \mathcal{O}(M^2)$ yield $\div \vv_{(0)}^{n+1} = \mathcal{O}(M^2)$ since $(\rho e + p)_{(0)}^n$ are constant.
    Consequently, for well-prepared initial data, the RS-IMEX FV scheme gives a consistent approximation of the limit equations as the Mach number tends to 0 independently of $\Delta t$.
\end{corollary}

\begin{proof}
    The proof can be done following the lines of the proof of Theorem \ref{theo:AP}.
    Due to $\rho^n$ and $\rho^{n+1}$ being a convex combination of the states $\rho_{1,\RS}$ and $\rho_{2,\RS}$, they are constant.
    As a consequence, we obtain \eqref{eq.corol.1}.
    Further, $\div \vv_{(0)}^{n+1} = \mathcal{O}(M^2)$ since $(\rho e + p)_{(0)}^n$ is constant.
    As the initial data are well-prepared, i.e. $\q^0 \in \Omega_{wp}^M$, we also obtain recursively $\q^n \in \Omega_{wp}^M$ for all successive time iterations $n > 0$.
\end{proof}

\section{Numerical results}\label{sec:NumRes}
In this section, we illustrate by numerical experiments theoretical properties of the first and second order RS-IMEX FV scheme, denoted respectively by RS-IMEX1 and RS-IMEX2, proposed in Section \ref{sec:NumScheme}.
All test cases are performed under the material CFL condition \eqref{eqn.CFL.cond} based on the eigenvalues of the explicit subsystem \eqref{eq.subsys.ex} which are of the order of the advection scale.
The initial conditions, if not mentioned otherwise, are given in dimensional form using the transformations \eqref{eq.scaled.quantities}, \eqref{eq.scaling_T} and the definition of the Mach numbers \eqref{eq.Mach}.
Whenever possible, we compare the numerical results obtained with our RS-IMEX FV schemes with an exact or explicit reference solution of the two-fluid model with single temperature \eqref{eq.SHTC.div}.
    Note that for a fully explicit scheme, the CFL condition depends directly on the Mach numbers $M_1,M_2$ and on the relaxation parameters $\tau^{(\vf)}, \tau^{(\ww)}$.
    Therefore, a fully explicit scheme which treats the relaxation source terms explicitly, is only comparable to the CFL condition \eqref{eqn.CFL.cond} of the RS-IMEX FV scheme for sonic and supersonic flows $M_1,M_2 \geq 1$ and slow relaxation processes $\tau^{(\vf)}, \tau^{(\ww)} \gg 1$.
    Thus, for test cases with well-prepared initial data, $\tau^{(\ww)} = \mathcal{O}(M^2)$, the RS-IMEX FV scheme allows significantly larger time steps in low Mach number regimes.
    Moreover, the case of pressure equilibrium which can be interpreted as “instantaneous” relaxation, is impossible to resolve with a purely explicit scheme. Thus, also in the case of
    an explicit scheme a nonlinear system has to be solved for $\vf_1$ to guarantee $p_1 = p_2$.
    This underlines the necessity of a semi-implicit AP scheme in such situations.

\subsection{Numerical convergence study}
\label{sec.convergence}
To verify the experimental order of convergence (EOC) we construct an exact solution of the
homogeneous two-fluid model single temperature system \eqref{eq.SHTC.div} given by a stationary
vortex.
It is obtained by considering zero radial velocities and a constant solution to be in angular direction, i.e.
\begin{equation}
    \label{eq.Vortex.assumptions}
    \ur = 0, \quad \wr = 0, \quad \frac{\pd}{\pd t} \left(\cdot\right)= 0, \quad \frac{\pd}{\pd \theta}\left(\cdot\right) = 0.
\end{equation}
In Appendix \ref{App:Polar}, the two-phase model \eqref{eq.SHTC.div} without the relaxation source terms is written in polar coordinates \eqref{eq:TFM_polar}.
Applying \eqref{eq.Vortex.assumptions}, it reduces to
\begin{subequations}
    \label{eq:TFM_vortex}
    \begin{eqnarray}
        &&\frac{\pd p}{\pd r} = \frac{\alpha\indd{1} \rho\indd{1} \uphi\indd{1}^2 + \alpha\indd{2} \rho\indd{2} \uphi\indd{2}^2}{r},\label{eq:TFM_vortex.p}\\
        && \frac{\pd}{\pd r} \left(\frac{\uphi\indd{1}^2 - \uphi\indd{2}^2}{2} + \mu_1 -
        \mu_2\right) - \uphi \left(\frac{1}{r} \frac{\pd}{\pd r}\left(r w_\theta\right)\right) = 0,
        \label{eq:TFM_vortex.w}
    \end{eqnarray}
\end{subequations}
with $p = \vf\indd{1} p\indd{1} + \alpha\indd{2} p\indd{2}, ~w_\theta = \uphi\indd{1} - \uphi\indd{2}, ~v_\theta = \mf\indd{1} \uphi\indd{1} + \mf\indd{2} \uphi\indd{2}$.
We set the phase velocities and the profile for the volume fraction to
\begin{equation*}
    \uphi\indd{l} = r v\indd{c,l} \exp\left(\nu\indd{\vv,l} (1- r^2)\right) \quad \text{and} \quad
    \vf\indd{1} = c_\alpha + \alpha_c \exp\left(\nu_\alpha (1-r^2)\right),
\end{equation*}
respectively.
This yields two equations for three unknowns $\rho\indd{1}, ~\rho\indd{2}$ and $T$.
To eliminate one unknown, we set $\rho\indd{2} = c\indd{\rho} \rho\indd{1}$ with $c\indd{\rho}$ being constant.
Then, the unknowns $\rho\indd{1}$ and $T$ can be determined via the following system of ordinary differential equations
\begin{equation}
    \label{eq.Vortex.ode}
    \begin{pmatrix}
        \displaystyle \vf\indd{1} \frac{\pd p\indd{1}}{\pd \rho\indd{1}} + \alpha\indd{2} c\indd{\rho} \frac{\pd p\indd{2}}{\pd \rho\indd{2}} &\displaystyle \vf\indd{1} \frac{\pd p\indd{1}}{\pd T} + \alpha\indd{2} \frac{\pd p\indd{2}}{\pd T} \\[3mm]
        \displaystyle\frac{\pd \mu\indd{1}}{\pd \rho\indd{1}} - c\indd{\rho}\frac{\pd \mu\indd{2}}{\pd \rho\indd{2}}&\displaystyle  \frac{\pd \mu\indd{1}}{\pd T} - \frac{\pd \mu\indd{2}}{\pd T} \\[3mm]
    \end{pmatrix}
    \begin{pmatrix}
        \displaystyle\frac{\pd \rho\indd{1}}{\pd r} \\[3mm]
        \displaystyle\frac{\pd T}{\pd r} \\[3mm]
    \end{pmatrix}
    =
    \begin{pmatrix}
        \displaystyle\frac{\alpha\indd{1} \rho\indd{1} \uphi\indd{1}^2 + \alpha\indd{2} \rho\indd{2} \uphi\indd{2}^2}{r} - p\indd{1} \frac{\pd \vf\indd{1}}{\pd r} + p\indd{2} \frac{\pd \vf\indd{1}}{\pd r} \\[3mm]
        \displaystyle -\frac{\pd}{\pd r} \left(\frac{\uphi\indd{1}^2 - \uphi\indd{2}^2}{2} \right) + \uphi \left(\frac{1}{r} \frac{\pd}{\pd r}\left(r w_\theta\right)\right)
    \end{pmatrix}.
\end{equation}
Applying the ideal gas law yields
\begin{subequations}
    \begin{eqnarray}
        &&\frac{\pd p\indd{1}}{\pd \rho\indd{1}} = (\gamma\indd{1} - 1)c\indd{v,1} T, \quad \frac{\pd p\indd{2}}{\pd \rho\indd{2}} = (\gamma\indd{2} - 1)c\indd{v,2} T\\[2mm]
        &&\frac{\pd p\indd{1}}{\pd T}= (\gamma\indd{1} - 1)c\indd{v,1} \rho\indd{1}, \quad \frac{\pd p\indd{2}}{\pd T} = (\gamma\indd{2} - 1)c\indd{v,2} \rho\indd{2} \\[2mm]
        &&\frac{\pd \mu\indd{1}}{\pd \rho\indd{1}}=\frac{(\gamma\indd{1}-1) c\indd{v,1}  T}{\rho\indd{1}}, \quad \frac{\pd \mu\indd{2}}{\pd \rho\indd{1}}= \frac{(\gamma\indd{2}-1) c\indd{v,2}  T}{\rho\indd{2}}\\[2mm]
        &&\frac{\pd \mu\indd{1}}{\pd T}= (\gamma\indd{1} -1) c\indd{v,1} -s\indd{1}, \quad \frac{\pd \mu\indd{1}}{\pd T}= (\gamma\indd{2} -1). c\indd{v,2} -s\indd{2}
    \end{eqnarray}
\end{subequations}
To obtain the initial condition on the computational domain $[-1,1]^2$, we integrate \eqref{eq.Vortex.ode} numerically with RK4, starting with the initial data $\rho\indd{1} = 1, \rho\indd{2} = 1, ~ T = 2$.
The parameters in order to obtain the velocities $v\indd{l,\theta}$ and the volume fraction $\alpha\indd{1}$
are set as
\begin{equation}
     c_\alpha = 0.4, \quad \alpha_c = 10^{-4}, \quad \nu_\alpha = 10, \quad v\indd{c,1} = 2 \cdot 10^{-5}, \quad v\indd{c,2} = 2.5 \cdot 10^{-5}, \quad \nu\indd{\vv,1} = 15, \quad \nu\indd{\vv,2} = 14.
\end{equation}
This setting yields two different phase velocities and consequently a non-zero relative velocity.
To obtain a vortex in the compressible flow regime, we assign the following material parameters
\begin{eqnarray}
    \label{eq.Vortex.param.1}
   \gamma\indd{1} = \frac{7}{5}, \quad \gamma\indd{2} = \frac{5}{3}, \quad c\indd{v,1} = 1, \quad c\indd{v,2} = 1.
\end{eqnarray}
The maximal Mach number for phase 1 is 0.62, for phase 2 it is 0.21 and the maximal mixture Mach number is 0.54. Consequently, the flow is compressible.

Since the sound speeds depend on the magnitude of the pressures which itself depend on $c_{v,l}$,
we scale $c_{v,l}$ with one over the Mach number $M$ to achieve flows in a desired Mach number regime.
In the next test case, we set $c_{v,l}/M^2$ which yields a maximum Mach
number of $0.018$ for the phase 1 and $0.014$ for the phase 2 and the mixture Mach number of $0.016$.
Setting
\begin{eqnarray}
    \label{eq.Vortex.param.2}
    \gamma\indd{1} = 2, \quad \gamma\indd{2} = 2.8, \quad c\indd{v,1} = 20, \quad c\indd{v,2} = 20,
\end{eqnarray}
the vortex flow is now weakly compressible.
Note that, according to Definition \ref{def.wp}, the initial data are ill-prepared since the phase densities are not constant.
However, we see from Tables \ref{tab.Vortex1} and \ref{tab.Vortex2} that the numerical scheme RS-IMEX2 FV converges with the expected EOC of two for both Mach number regimes.
The results are obtained with a material CFL condition \eqref{eqn.CFL.cond} with $\nu = 0.25$.


\begin{table}
    \centering
    \begin{tabular}{ccccccccc}
        & 16 & & 	32 & &	64 & &	128 &
        \\ \hline
        $\alpha\indd{1}$  & 6.31E-03 & ---& 	1.15E-03 & 2.45 &	2.20E-04	& 2.38 & 4.69E-05 &  2.23\\[2mm]
        $\rho\indd{1}	$   & 2.98E-02 & ---& 	9.44E-03 & 1.65 &	2.24E-03	& 2.07 & 3.94E-04 &  2.50\\[2mm]
        $\rho\indd{2}	$   & 2.78E-02 & ---& 	8.21E-03 & 1.75 &	1.60E-03	& 2.35 & 2.56E-04 &  2.64\\[2mm]
        $\vv\indd{1,1}$	& 5.51E-02 & ---& 	1.40E-02 & 1.98 &	2.44E-03	& 2.51 & 3.33E-04 &  2.87\\[2mm]
        $\vv\indd{2,1}$	& 5.51E-02 & ---& 	1.40E-02 & 1.97 &	2.45E-03	& 2.51 & 3.41E-04 &  2.84\\[2mm]
        $\vv\indd{1,2}$	& 6.85E-02 & ---& 	1.58E-02 & 2.11 &	2.50E-03	& 2.65 & 3.41E-04 &  2.87\\[2mm]
        $\vv\indd{2,2}$	& 6.85E-02 & ---& 	1.58E-02 & 2.11 &	2.49E-03	& 2.66 & 3.51E-04 &  2.82\\[2mm]
        $T$         & 4.45E-02 & ---& 	1.84E-02 & 1.27 &	3.92E-03	& 2.23 & 6.54E-04 &  2.58\\[1mm]\hline
    \end{tabular}

    \caption{Two-fluid stationary vortex: $L^1$ error and EOC for the second order RS-IMEX FV scheme in the compressible regime with parameters given in \eqref{eq.Vortex.param.1}.}
    \label{tab.Vortex1}
\end{table}

\begin{table}
    \centering
    \begin{tabular}{ccccccccc}
        & 16 & & 	32 & &	64 & &	128 &
        \\ \hline
        $\alpha\indd{1}$  & 6.53E-03 & ---& 	1.29E-03 & 2.34 &	2.57E-04	& 2.32 & 5.42E-05 &  2.24\\[2mm]
        $\rho\indd{1}	$   & 6.95E-03 & ---& 	1.70E-03 & 2.02 &	4.68E-04	& 1.86 & 1.19E-04 &  1.97\\[2mm]
        $\rho\indd{2}	$   & 2.73E-02 & ---& 	1.20E-03 & 4.50 &	3.23E-04	& 1.89 & 8.17E-05 &  1.98\\[2mm]
        $\vv\indd{1,1}$	& 4.11E-01 & ---& 	1.32E-02 & 4.95 &	2.09E-03	& 2.66 & 3.08E-04 &  2.76\\[2mm]
        $\vv\indd{2,1}$	& 4.11E-01 & ---& 	1.33E-02 & 4.95 &	2.09E-03	& 2.66 & 3.08E-04 &  2.76\\[2mm]
        $\vv\indd{1,2}$	& 4.17E-01 & ---& 	2.64E-02 & 3.98 &	6.14E-03	& 2.10 & 1.41E-03 &  2.12\\[2mm]
        $\vv\indd{2,2}$	& 4.17E-01 & ---& 	2.64E-02 & 3.98 &	6.14E-03	& 2.10 & 1.41E-03 &  2.12\\[2mm]
        $T$         & 2.50E-02 & ---& 	3.20E-03 & 2.96 &	8.24E-04	& 1.95 & 2.01E-04 &  2.03\\[1mm]\hline
    \end{tabular}

    \caption{Two-fluid stationary vortex: $L^1$ error and EOC for the second order RS-IMEX FV scheme in the
    weakly-compressible regime with parameters given in \eqref{eq.Vortex.param.2}.}
    \label{tab.Vortex2}
\end{table}

\subsection{1D Riemann Problems}
\label{sec.rp}

\begin{figure}
    \vspace{-15mm}
    \centering
    \includegraphics[scale=0.38]{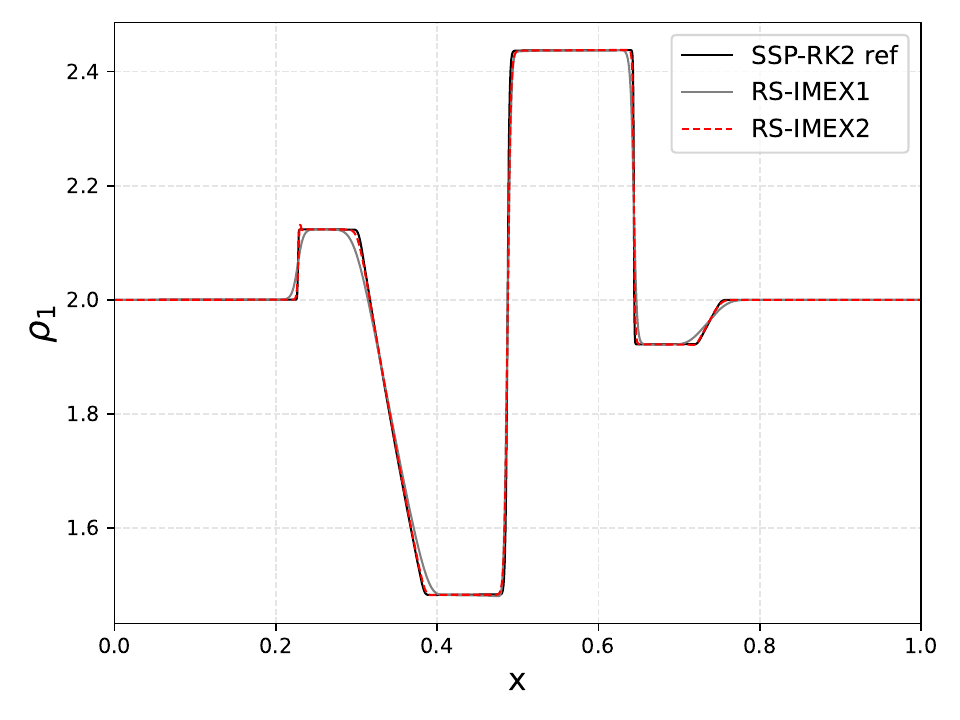} \includegraphics[scale=0.38]{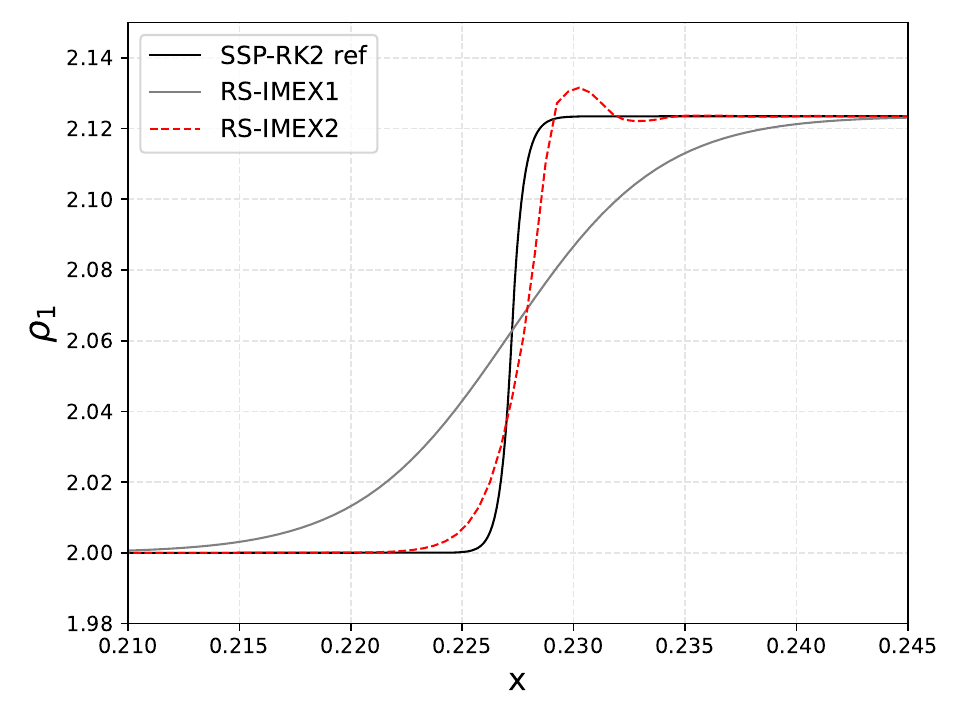}\\
    \includegraphics[scale=0.38]{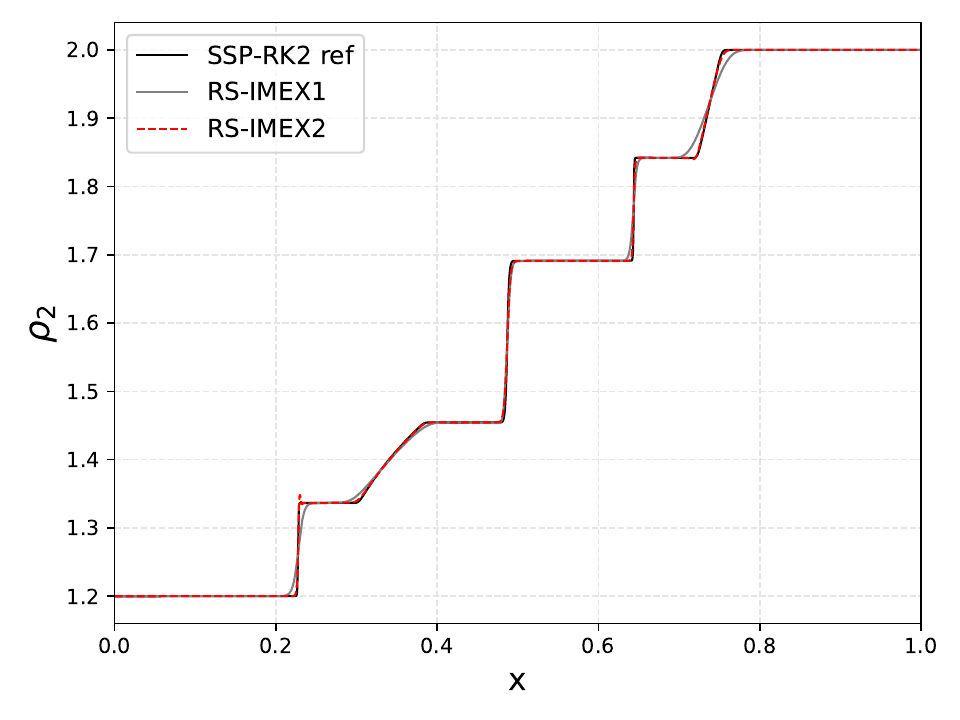}
    \includegraphics[scale=0.38]{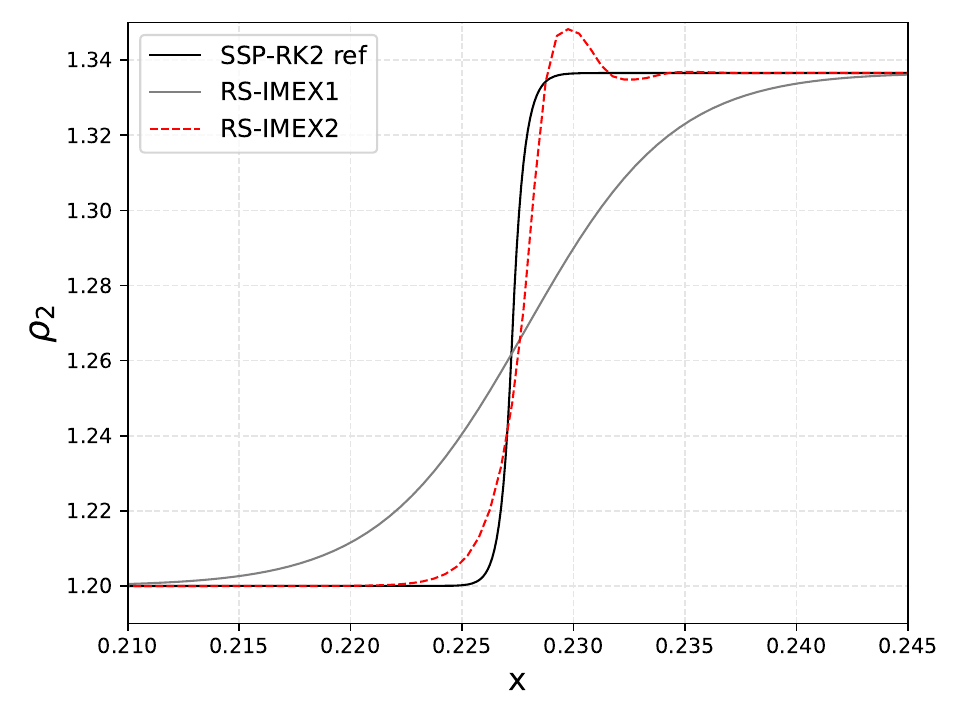}\\
    \includegraphics[scale=0.38]{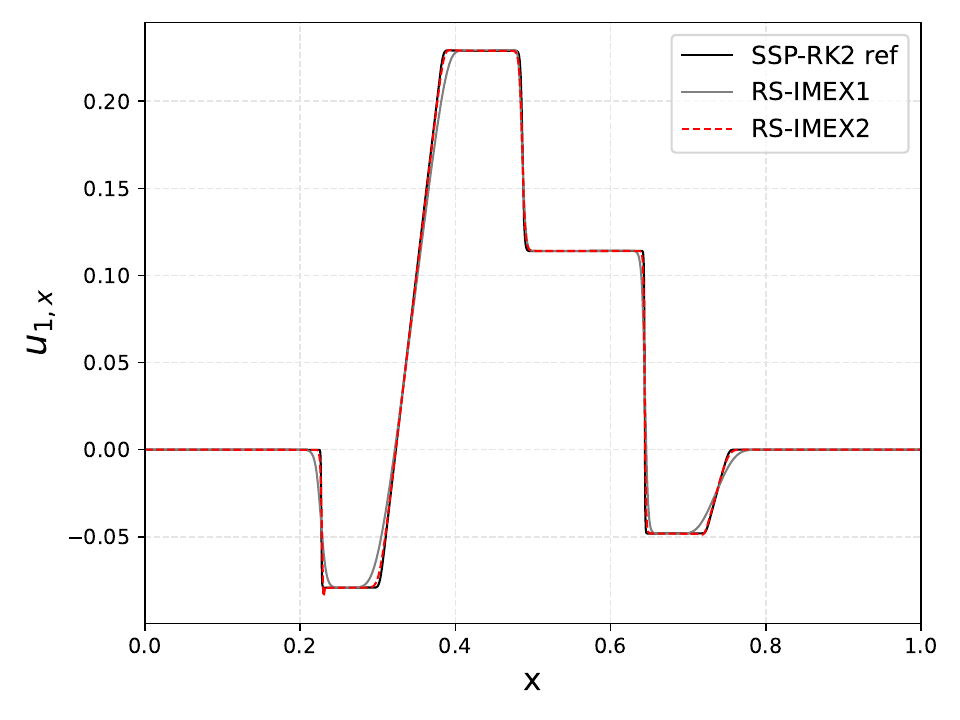}
    \includegraphics[scale=0.38]{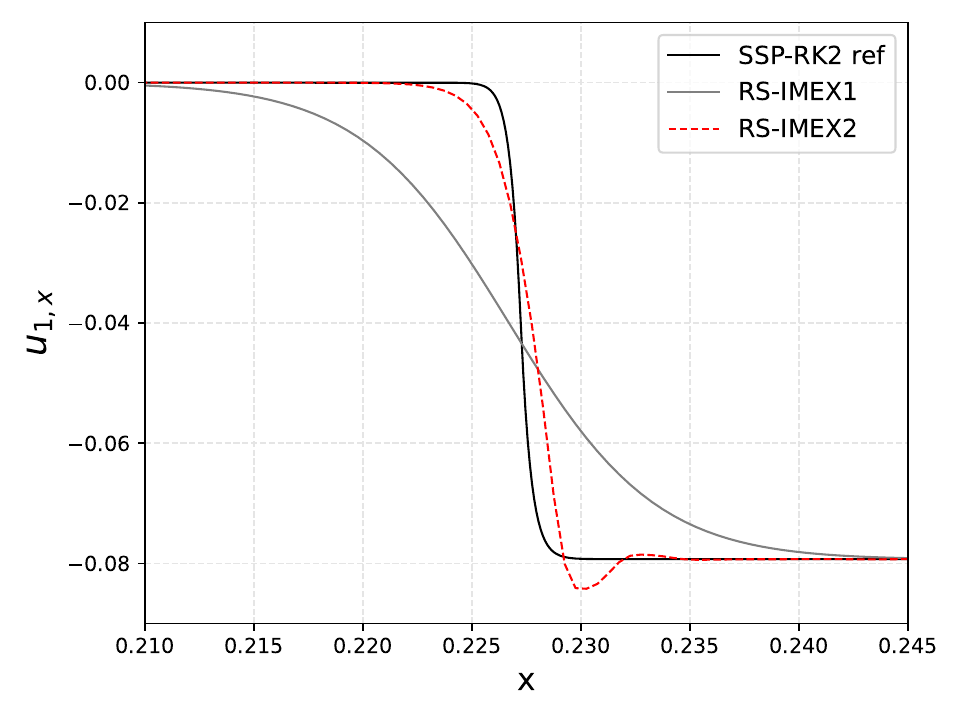}\\
    \includegraphics[scale=0.38]{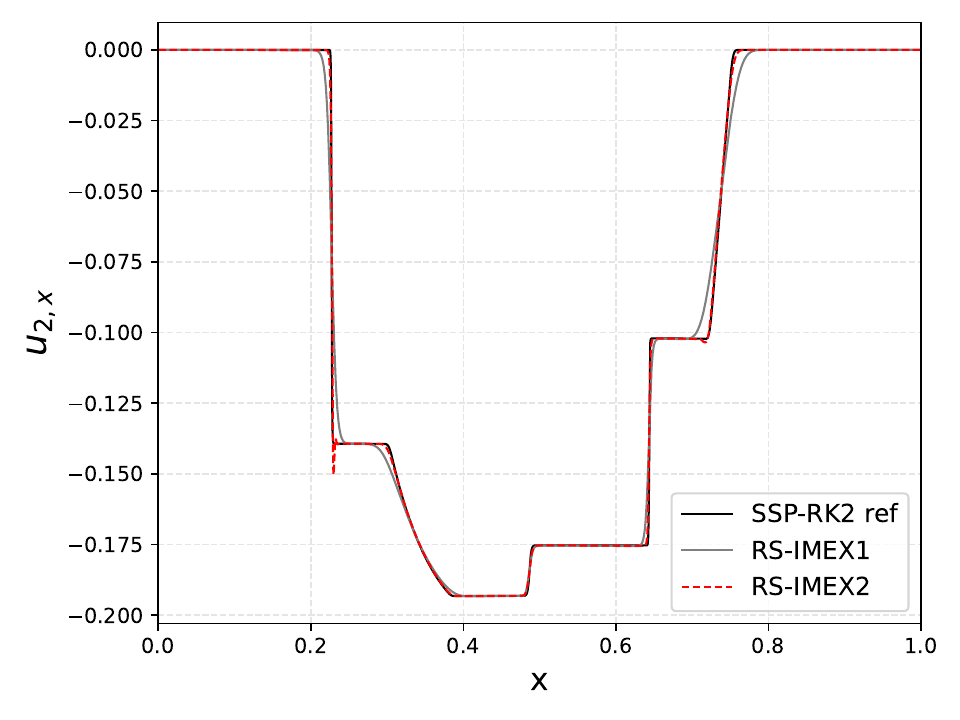}
    \includegraphics[scale=0.38]{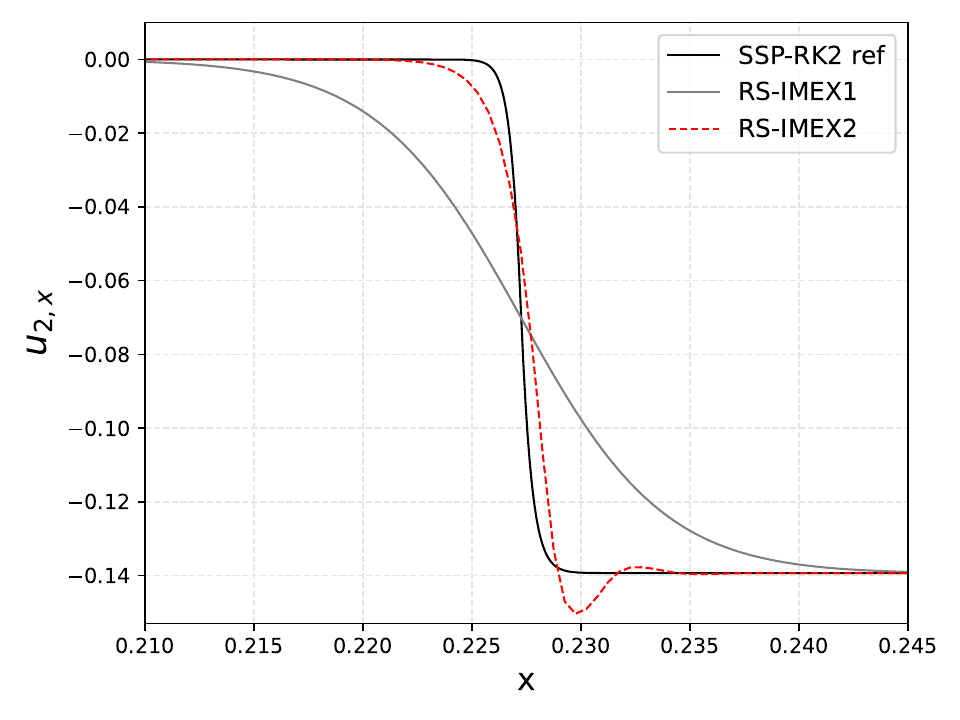}\\
    \includegraphics[scale=0.38]{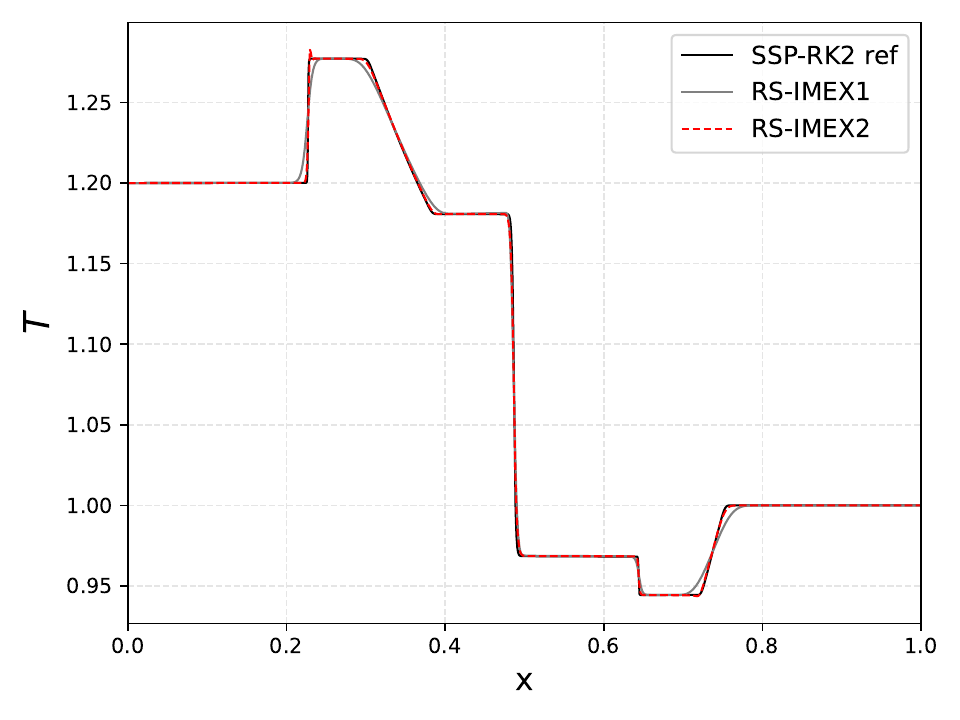}
    \includegraphics[scale=0.38]{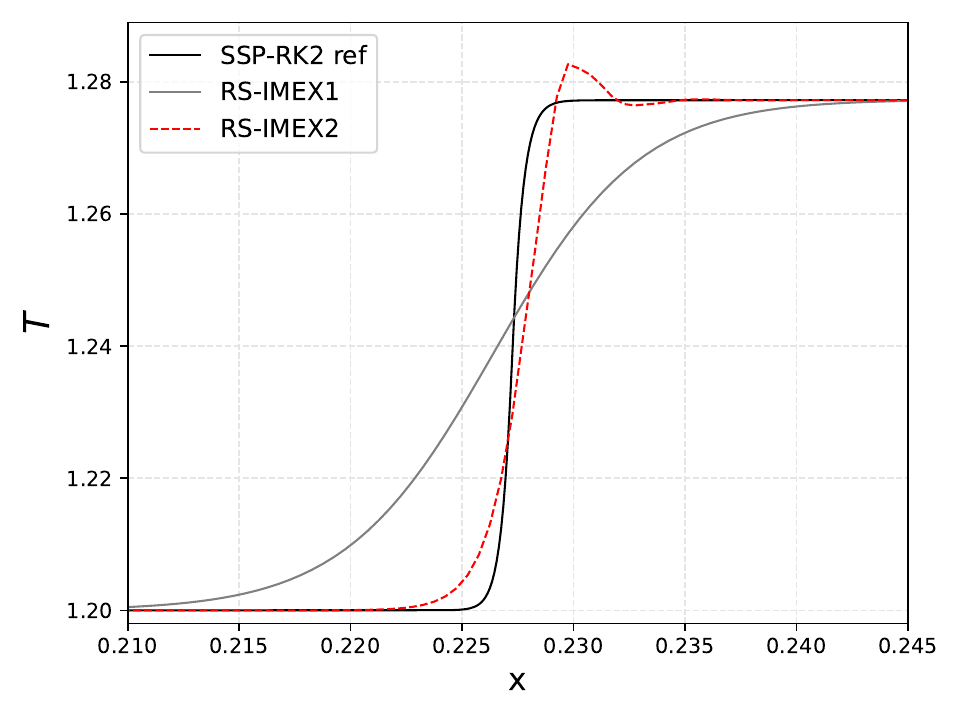}
    \caption{Numerical solutions of the homogeneous Riemann problem RP1 obtained at time $T_f=0.2$ with constant volume fraction without relaxation source terms using the new first and second order RS-IMEX FV schemes.
    The reference solution is computed by the explicit second order SSP-RK2 FV scheme. From top
    left to bottom right: Phase densities $\rho_1,~\rho_2$, phase velocities $\vv_{1,1},~\vv_{2,1}$
    and temperature $T$. Left: Computational domain $x\in[0,1]$. Right: Zoom on
    $x\in[0.21,0.245]$.}
\label{fig.rp1}
\end{figure}

\begin{table}[t!]
    \centering
    \renewcommand{\arraystretch}{1.25}
    \begin{tabular}{lclcccccc}
        Test                 & $T_f$ &   state & $\vf$ & $\rho_{1}$ & $\rho_2$ & $ \vv_{1,1} $ & $ \vv_{2,1}$ & $T$  \\\hline \hline
        \multirow{2}{*}{RP1} & \multirow{2}{*}{0.2}& left  & 0.3 & 2  & 1.2      & 0         & 0        & 1.2     \\
        &              & right & 0.3 & 2  & 2        & 0         & 0        & 1     \\\hline
        \multirow{2}{*}{RP2} & \multirow{2}{*}{0.2}    & left  & 0.7 & 1  & 2        & -1        & -1       & 1     \\
        &                   & right & 0.3 & 1  & 2        &  1        &  1       & 1     \\\hline
    \end{tabular}
    \caption{Initial condition for the 1D Riemann problems presented in Section \ref{sec.rp} with $\gamma\indd{1} = 1.4$ and $\gamma\indd{2} = 2$ and $c_{v,1} = c_{v,2} = 1$ on the domain $[0,1]$ with initial jump at $x = 0.5$. }
    \label{tab.initRP}
\end{table}
To test the first and second order versions of the RS-IMEX FV scheme in a high Mach number regime, we
consider two Riemann Problems (RPs) for the homogeneous system \eqref{eq.SHTC.div}
omitting the pressure relaxation source term acting on the volume fraction \eqref{eqn.alphaMS.div}
and the friction source term in the relative velocity equation \eqref{eqn.relvelMS.div}.
The initial configuration on the domain $[0,1]$ is given in Table \ref{tab.initRP} and the initial jump position is located at $x=0.5$.
The first RP (RP1) consists of an initial jump in density of phase two and the temperature where the volume fraction is kept constant.
The second RP (RP2) is a double rarefaction test with an initial jump in the volume fraction
resulting in a discontinuous mixture density and internal energy.
In Figure \ref{fig.rp1}, we compare the results for the first and second order RS-IMEX FV schemes using 2000 cells and the material CFL condition \eqref{eqn.CFL.cond}.
For the RS-IMEX1 FV scheme we set $\nu = 0.8$ and for the RS-IMEX2 FV scheme $\nu = 0.4$.
This results in $\Delta t = 4 \cdot 10^{-4}$ and $ 2 \cdot 10^{-4}$, respectively.
The reference solution was computed by a second order explicit SSP-RK2 FV scheme using 10000 cells resulting in $\Delta t = 7.65 \cdot 10^{-6}$.
Note that the CFL condition for the explicit scheme is dictated by the fastest wave speed arising in the model which depends on the sound speeds of the respective phases, \cite{RomTor2004}.
Moreover a comparable time step of the explicit scheme for 2000 cells is $3.3\cdot 10^{-5}$ which is 10 times smaller than the one used for the IMEX schemes.
Since there are no shear processes present, RP1 consists of 5 waves.
The wave speeds and positions produced by both RS-IMEX FV schemes are in good agreement with the reference solution, where the first order scheme is more diffusive on the fast waves than the second order scheme, for which we can observe small oscillations on the outermost fast travelling waves, see right panel of Figure \ref{fig.rp1}.
Their appearance is local and does not impair the results on the material wave which is the focus of the simulation and is captured accurately by both RS-IMEX FV schemes since the chosen time step is oriented towards its accurate capturing only.
Moreover, the phenomenon of spurious oscillations around discontinuities is a known problem for higher order numerical schemes.
To fully resolve all waves, an acoustic time step can be chosen or additional artificial viscosity can be added in the explicit upwind part at the cost of more diffusive material waves.

\begin{figure}
    \vspace{-15mm}
    \centering
    \includegraphics[scale=0.38]{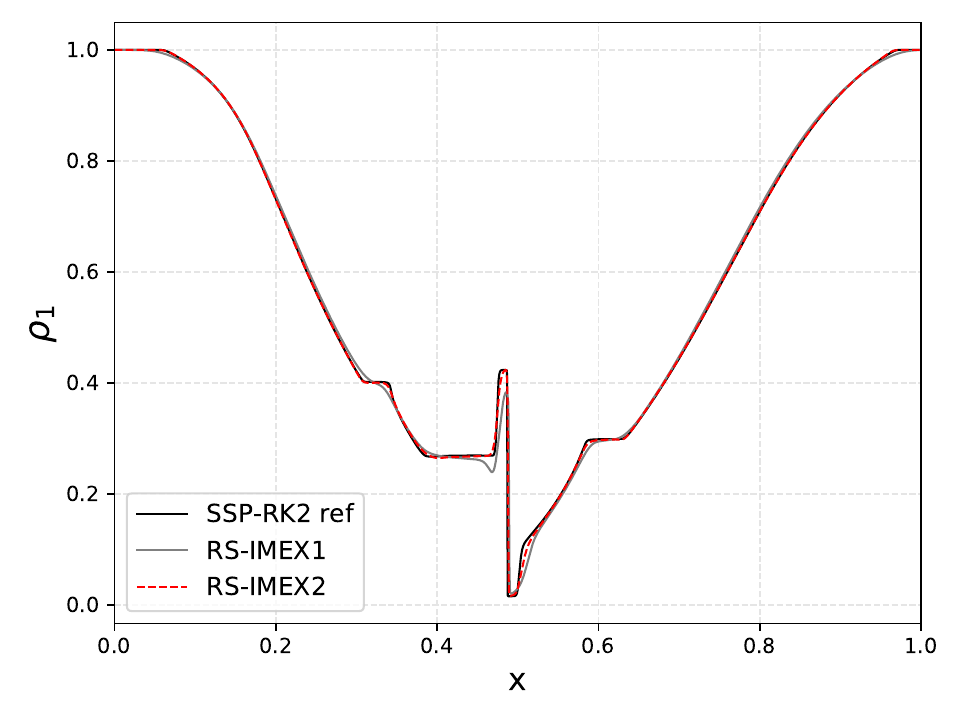} \includegraphics[scale=0.38]{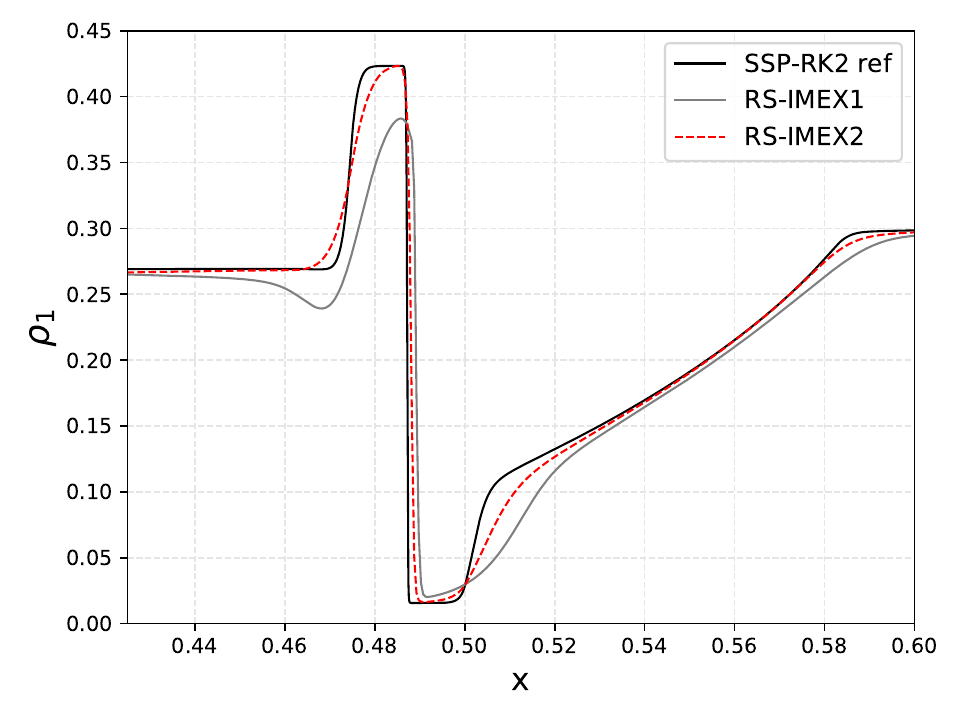}\\
    \includegraphics[scale=0.38]{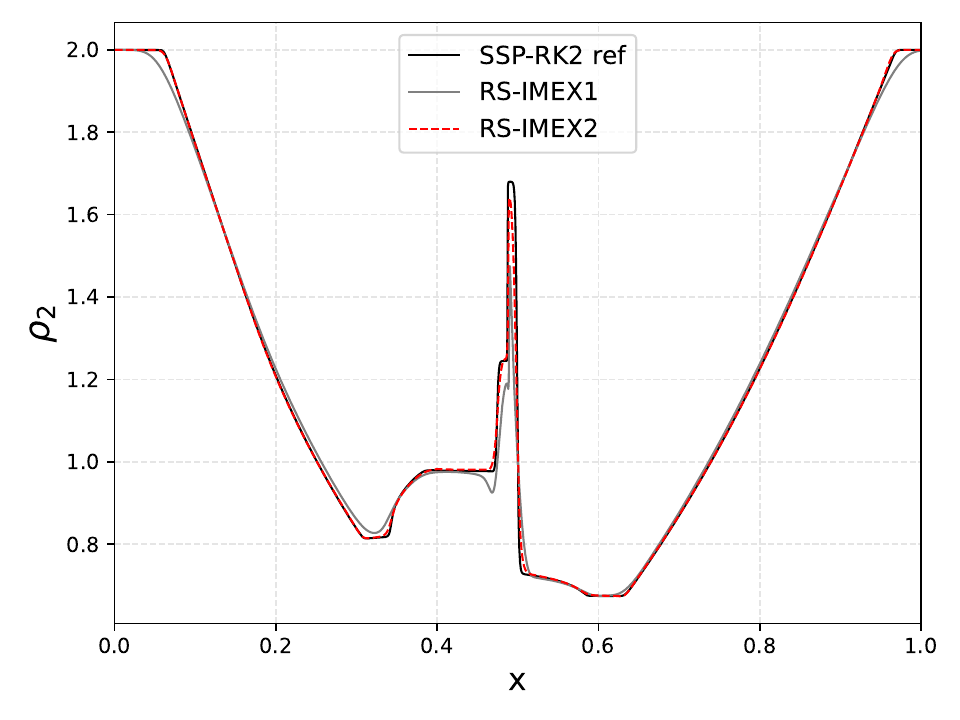}
    \includegraphics[scale=0.38]{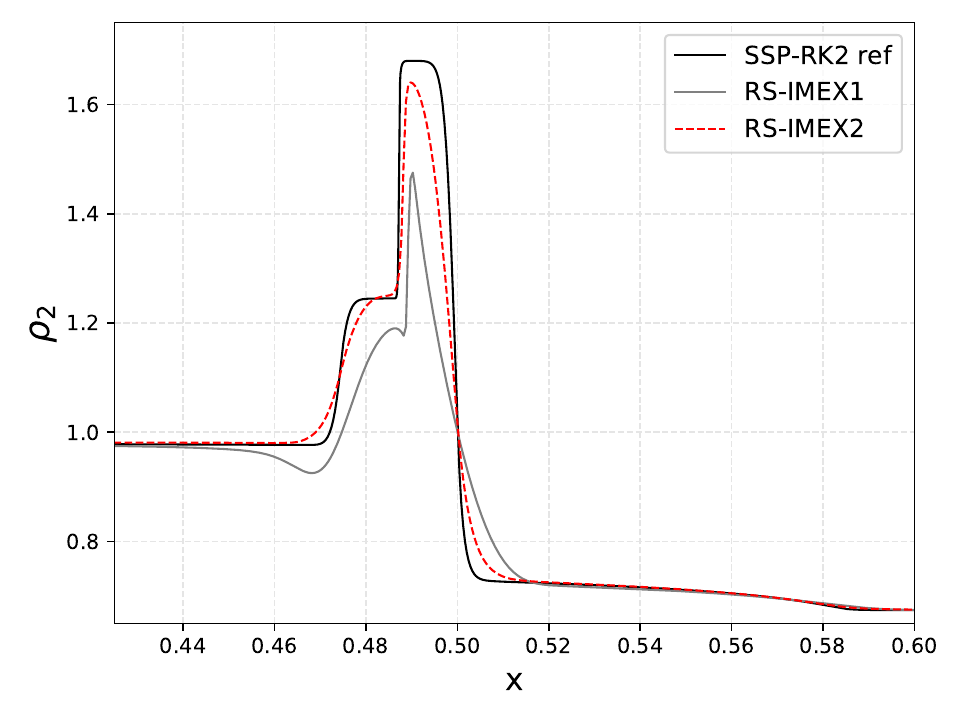}\\
    \includegraphics[scale=0.38]{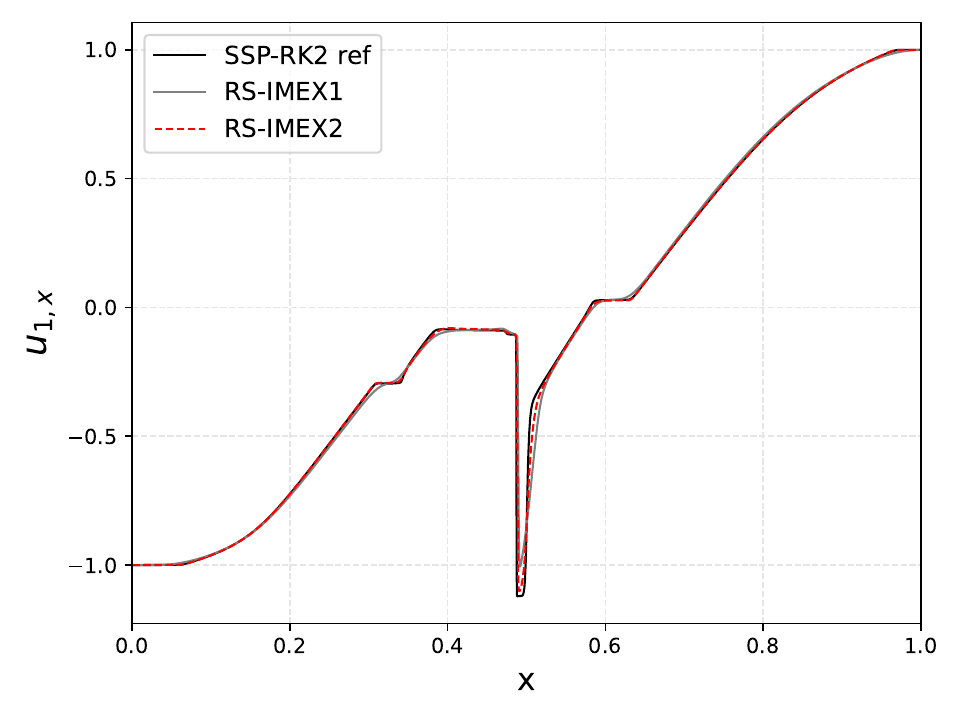}
    \includegraphics[scale=0.38]{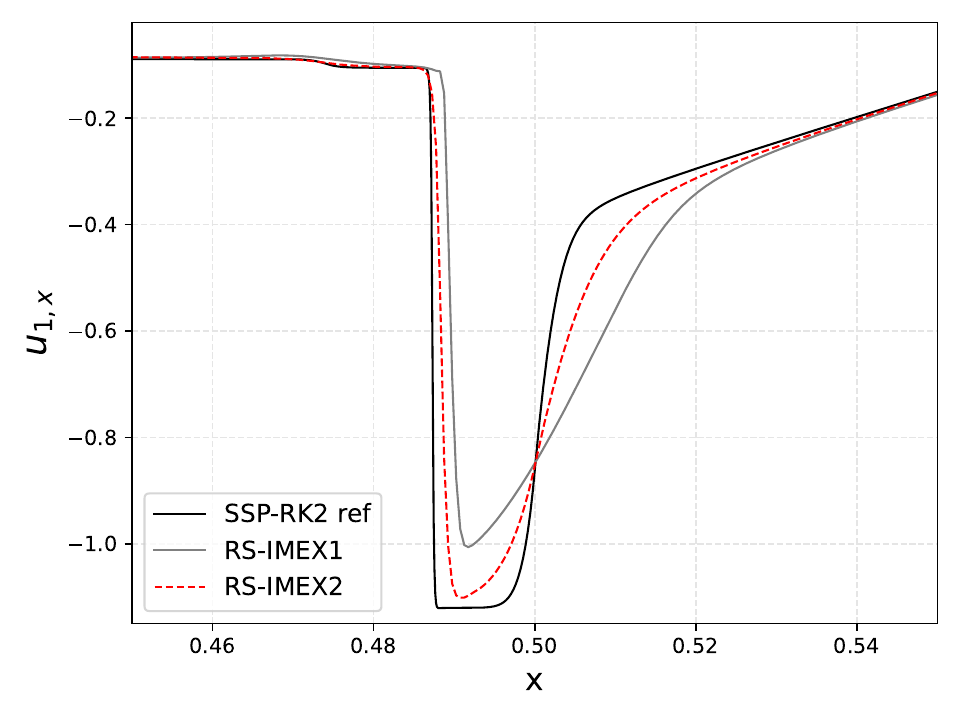}\\
    \includegraphics[scale=0.38]{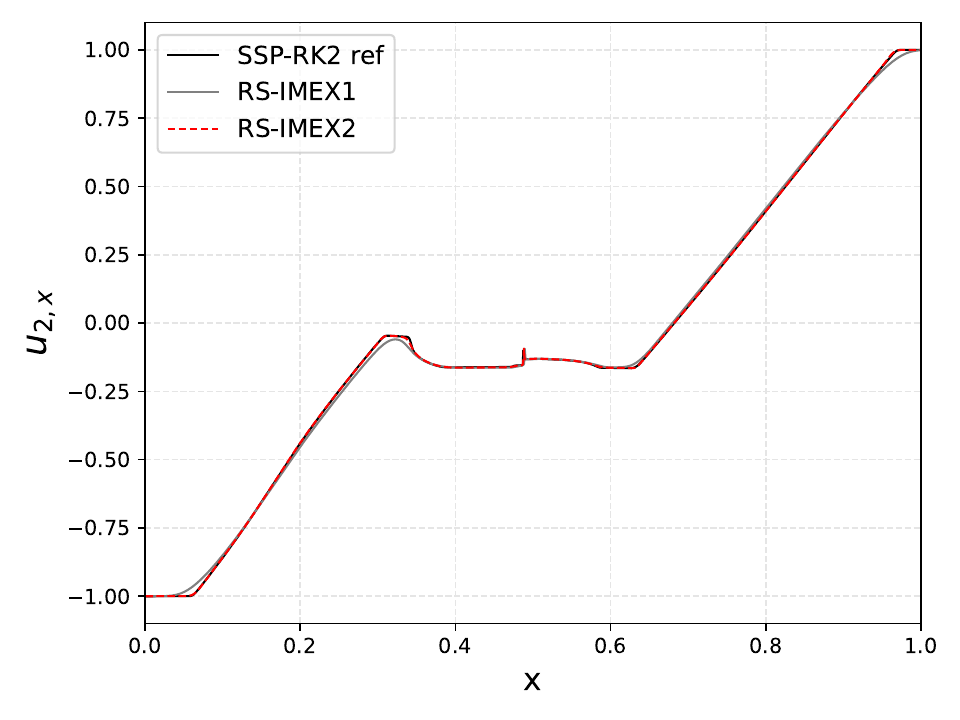}
    \includegraphics[scale=0.38]{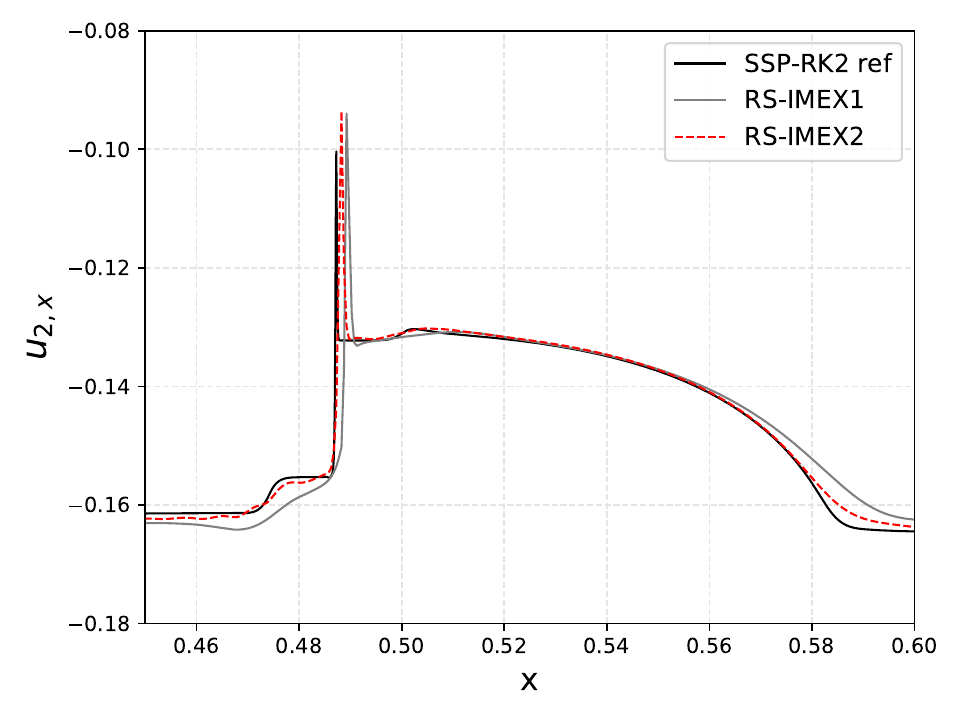}\\
    \includegraphics[scale=0.38]{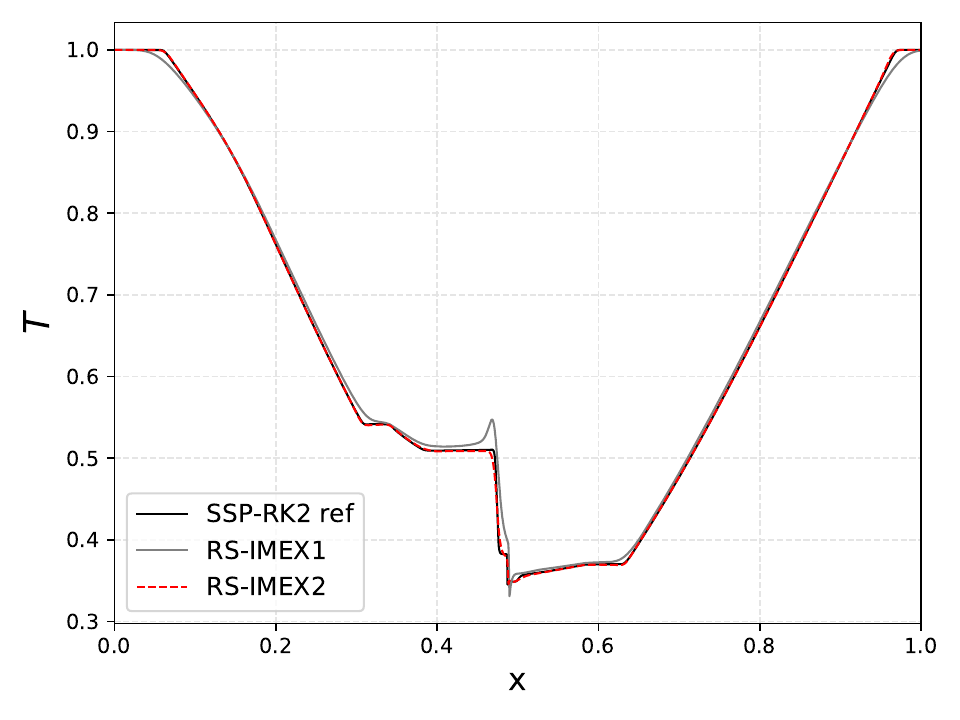}
    \includegraphics[scale=0.38]{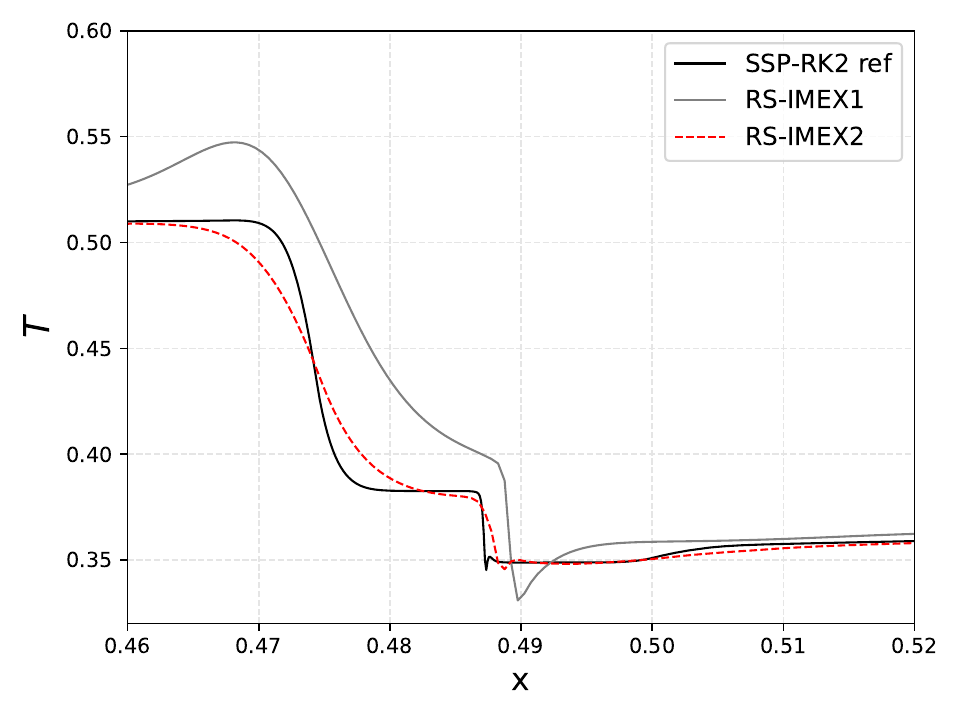}
    \caption{Numerical solutions of the homogeneous Riemann problem RP2 obtained at time $T_f=0.2$ with an initial jump in the volume fraction without relaxation source terms using the new first and second order RS-IMEX FV schemes.
    The reference solution is computed by the explicit second order SSP-RK2 FV scheme. From top to bottom: Phase densities $\rho_1,~\rho_2$, phase velocities $\vv_{1,1},~\vv_{2,1}$ and temperature $T$. Left: Computational domain $[0,1]$. Right: Zoom on material waves. }
    \label{fig.rp2}
\end{figure}

The wave structure of RP2 is more intricate, as can be seen in Figure \ref{fig.rp2}.
This is due to the initial jump in the volume fraction.
It consists of the contact wave associated with $\vf\indd{1}$ and three waves traveling to the left of the boundaries of the domain and two waves to the right.
Note that due to the single temperature assumption, the wave propagation is not symmetric.
We can observe, that the first order RS-IMEX1 FV scheme is too diffusive in order to capture the complicated sequence of slower waves.
On the other hand, the second order RS-IMEX2 FV scheme shows a great improvement in the capturing of the slower waves near the initial jump position.

\subsection{Advection of a Bubble}

\begin{figure}[ht!]
    \centering
    \includegraphics[scale=0.45]{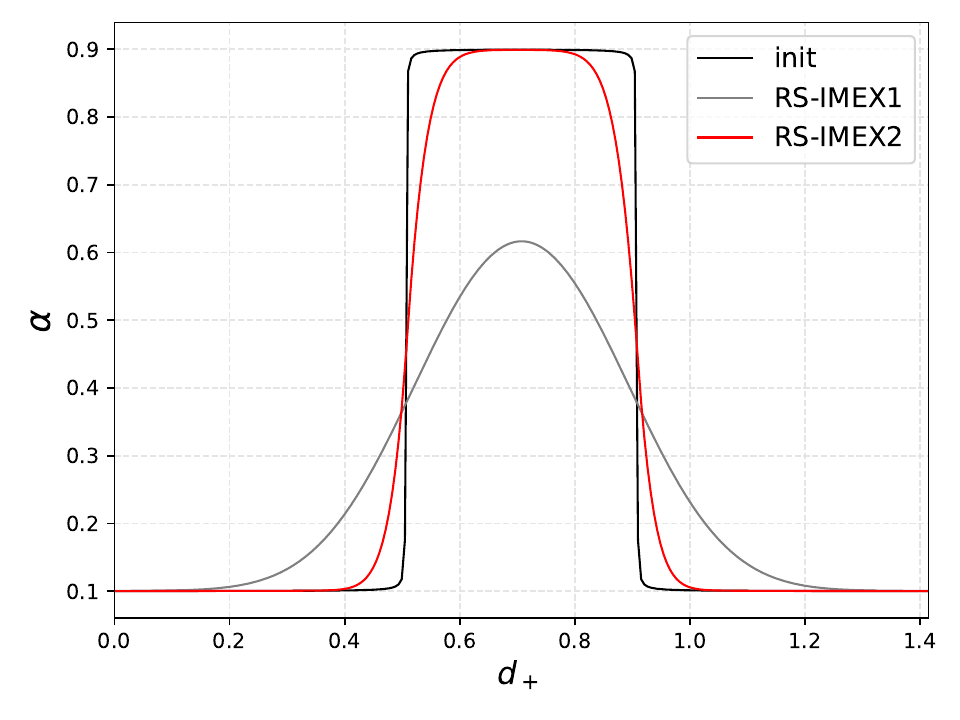}
    \includegraphics[scale=0.45]{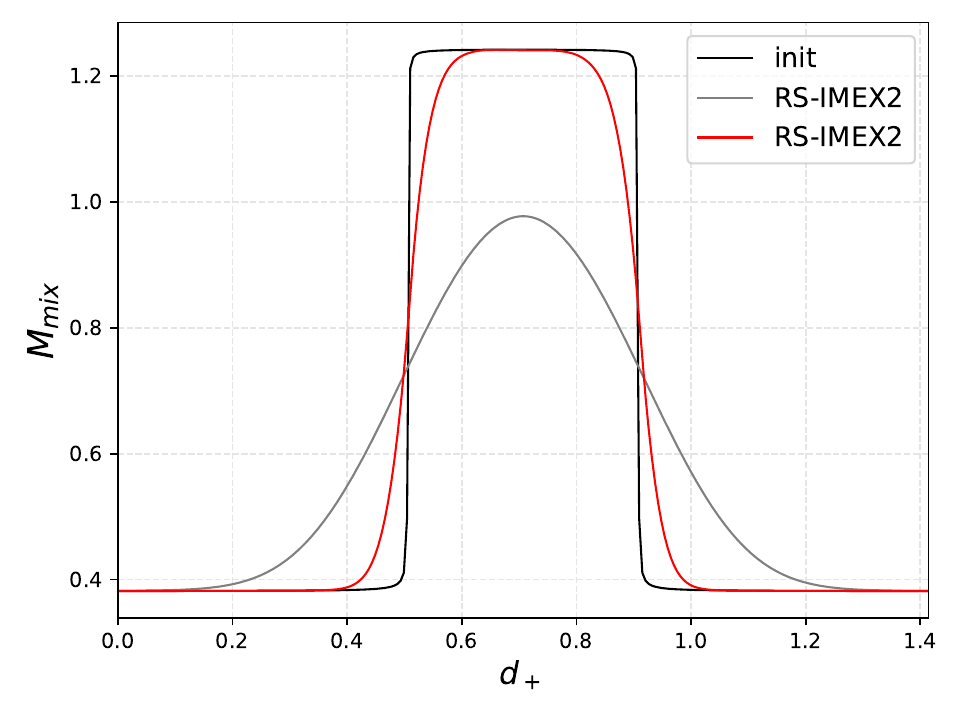}\\
    \includegraphics[scale=0.45]{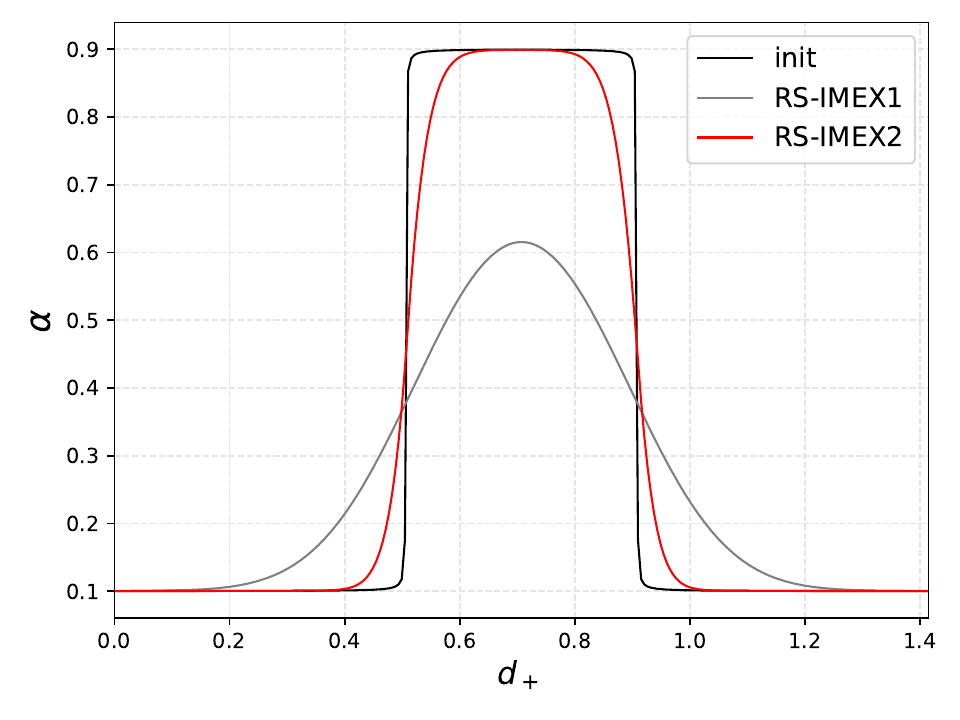}
    \includegraphics[scale=0.45]{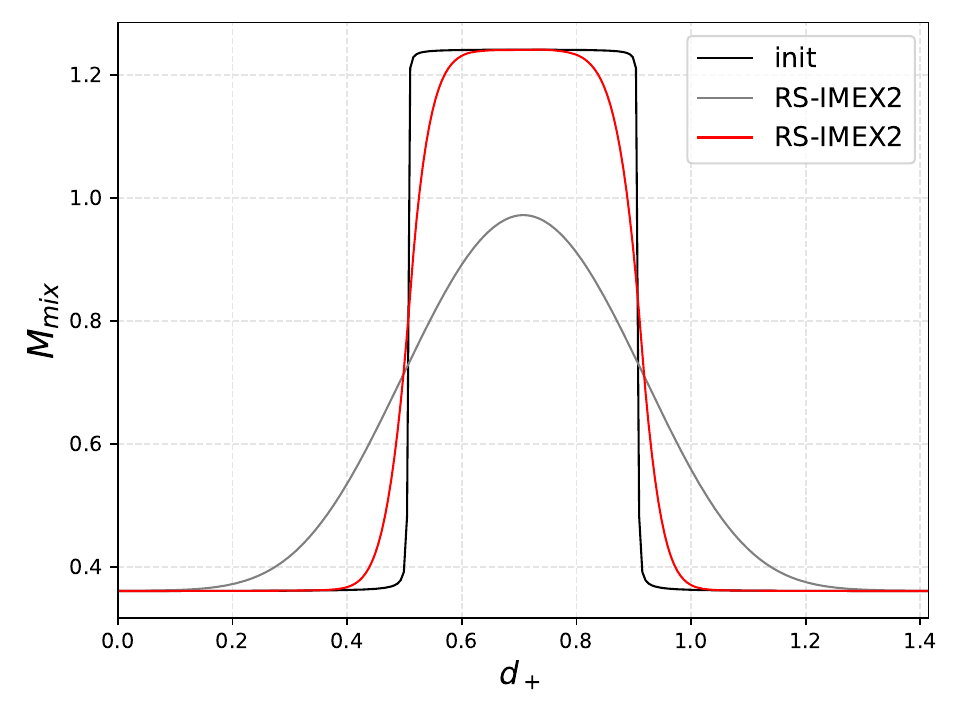}
    \caption{Numerical solutions of the diagonally advected bubble obtained at time $T_f=1$ by the first and second order RS-IMEX FV schemes displayed along the diagonal from $[0,0]$ to $[1,1]$. Top: Case $\mathcal{C} = 10$. Bottom: Case $\mathcal{C} = 50$. Left: Volume fraction $\vf\indd{1}$. Right: Mixture Mach number $M_{\text{mix}}$.}
    \label{fig:AdvecBubble}
\end{figure}

We consider a diagonally advected bubble initially centered at $(x_0,y_0) = (0.5,0.5)$ with the radius $r_0 = 0.2$.
The computational domain is set to $[0,1]\times[0,1]$ and is discretized by $256\times256$ rectangular mesh cells. Further we apply periodic boundary conditions.
The velocity fields are given by
\begin{equation}
    \vv_1 = (1,1)^T, \vv_2 = (1,1)^T.
\end{equation}
The bubble of phase 1 is moved through a second phase 2 which is modeled by a change in the volume fraction given in dependence of the radius $r = \sqrt{(x - x_0)^2 + (y - y_0)^2}$.
The initial volume fraction is given by
\begin{equation}
    \vf\indd{1}(r,0) = {(\vf_L - \vf_R)} \frac{\arctan\left(-\theta\left(r-r_0\right)\right)}{\pi} + \frac{(\vf_L + \vf_R)}{2},
\end{equation}
where $\vf_L = 0.9, ~\vf_R = 0.1$ and $\theta = 2000$. The parameter $\theta$ indicates the
diffusivity of the interface in the initial data.
To have a bubble, that is initially in pressure equilibrium, we set $\rho_2$ such that $p\indd{1} = p\indd{2}$, i.e.
\begin{equation}
    \label{eq.rho2_thermal_eq}
    \rho\indd{2} = \frac{(\gamma\indd{1}-1) c\indd{v,1}}{(\gamma\indd{2} - 1) c\indd{v,2}}\rho\indd{1},
\end{equation}
where $\rho\indd{1} = 2$.
To ensure that the phase-pressure equilibrium holds during the simulation, we set the relaxation
parameter $\tau^{(\vf)} = 10^{-16}$, i.e. ``instantaneous'' pressure relaxation.
Further, we set the initial temperature to $T=2$.
Finally, we set $\gamma\indd{1} = 1.4, ~\gamma\indd{2} = 2$, $c_{v,1} = 1$ and $c_{v,2}$ in accordance with \eqref{eq.Mach_C} by
\begin{equation}
    \label{eq:cv_Mach}
    c_{v,2} = \frac{\gamma_1 (\gamma_1 - 1) c_{v,1}}{\gamma_2 (\gamma_2 - 1)}\mathcal{C}^2.
\end{equation}
As before, $\mathcal{C}$ denotes the ratio between the Mach numbers and will be used to adjust the flow regimes.
Note, that the relative velocity equation for initially constant chemical potentials does not reduce to pure advection, but creates perturbations in $\ww$.
Therefore, to decrease these perturbations which can interfere with the advection of the bubble, we assume a high friction by setting $\tau^{(\ww)} = 10^{-8}$ and $\tau^{(\ww)} = 10^{-12}$ depending on the Mach number regimes associated with $\mathcal{C} = 10$ and $\mathcal{C}=50$, respectively.
This leads to the Mach numbers $M_{1,\max} = 1.336$ and $M_{2,\max} = 1.336\cdot10^{-1}$ for the first case and $M_{1,\max} = 1.336$ and $M_{2,\max} = 2.67\cdot10^{-2}$ for the second case.
The bubble is evolved up to the final time $T_f=1$ when the bubble is back in its initial position.
In Figure \ref{fig:AdvecBubble} the volume fraction $\vf$ together with the mixture Mach number \eqref{eq.Mach} is plotted along the diagonal for the first and second order schemes.
Both schemes use a material CFL condition \eqref{eqn.CFL.cond} with $\nu = 0.5$ for the first order scheme and $\nu = 0.25$ for the second order scheme.
The numerical solutions are in good agreement with the initial data.
The RS-IMEX1 FV scheme is quite diffusive whereas the RS-IMEX2 FV scheme captures well the initial configuration.
Note further, that the mixture Mach number changes rapidly from $\approx 1.23$ inside the bubble to $\approx 0.38$ outside of the bubble and from $\approx 1.24$ to $\approx 0.36$ for $\mathcal{C} = 10$ and $\mathcal{C}=50$, respectively.
Even though the phase Mach number $M_2$ is significantly smaller in the second case, the mixture Mach number does not change due to the averaging with respect to the mass fraction.
It is therefore not a good indicator to the individual Mach number regimes that determine the scales in the model.

\subsection{Kelvin Helmholtz instability}

We modify a set-up from \cite{MicRoepEde2015,McNalLyrPas2012} for the single phase Euler equations to the single temperature two-phase model \eqref{eq.SHTC.div}.
It describes two phases flowing in opposite directions which creates the Kelvin Helmholtz instability.
We apply periodic boundary conditions and set the computational domain to $[0,1]\times[0,1]$.
The two fluids are characterized by $\gamma_1 = 2$ and $\gamma_2 = 1.4$, respectively.
Further $\rho_1 = 1$ and $\rho_2$ is set according to \eqref{eq.rho2_thermal_eq} in such a way that the initial condition is in pressure equilibrium.
Furthermore, we require both fluids to have the same Mach number.
We set $c\indd{v,1} = 1/\varepsilon^2$ and $c\indd{v,2}$ with $\mathcal{C} = 1$ according to \eqref{eq:cv_Mach}.
Setting $T=12.5$ with $\varepsilon =1$ yields the maximal initial Mach number $M = 10^{-1}$, and choosing
$\varepsilon = 0.1$ yields the maximal initial Mach number $M = 3\cdot10^{-2}$.
To ensure that the flow stays in pressure equilibrium, we set $\tau^{(\vf)} = 10^{-16}$ and in accordance with the well-prepared initial data, we set $\tau^{(\ww)} = M^2$.
Initially, we choose the same phase velocities, defined as
\begin{equation}
    \vv_{1,1} = \vv_{2,1} =
    \begin{cases}
        v_L - v_m \exp((y-0.25)/L), &\quad \text{if} \quad 0 \leq y < 0.25 \\
        v_R + v_m \exp(-(y-0.25)/L), &\quad \text{if} \quad 0.25 \leq y < 0.5 \\
        v_R + v_m \exp((y-0.75)/L), &\quad \text{if} \quad 0.5 \leq y < 0.75 \\
        v_L - v_m \exp(-(y-0.75)/L), &\quad \text{if} \quad 0.75 \leq y \leq 1 \\
    \end{cases},
\end{equation}
where $v_L = 0.5$, $v_R = - 0.5$, $v_m = (v_L - v_R)/2$ and $L=0.025$.
In $y$-direction we apply an initial perturbation $\vv_{1,2} = \vv_{2,2} = 10^{-2}\sin(4 \pi x)$ which yields an initial relative velocity $\ww = 0$ and divergence free velocity field $\div \vv =0$.

\begin{figure}[t!]
    \centering
\includegraphics[scale=0.75]{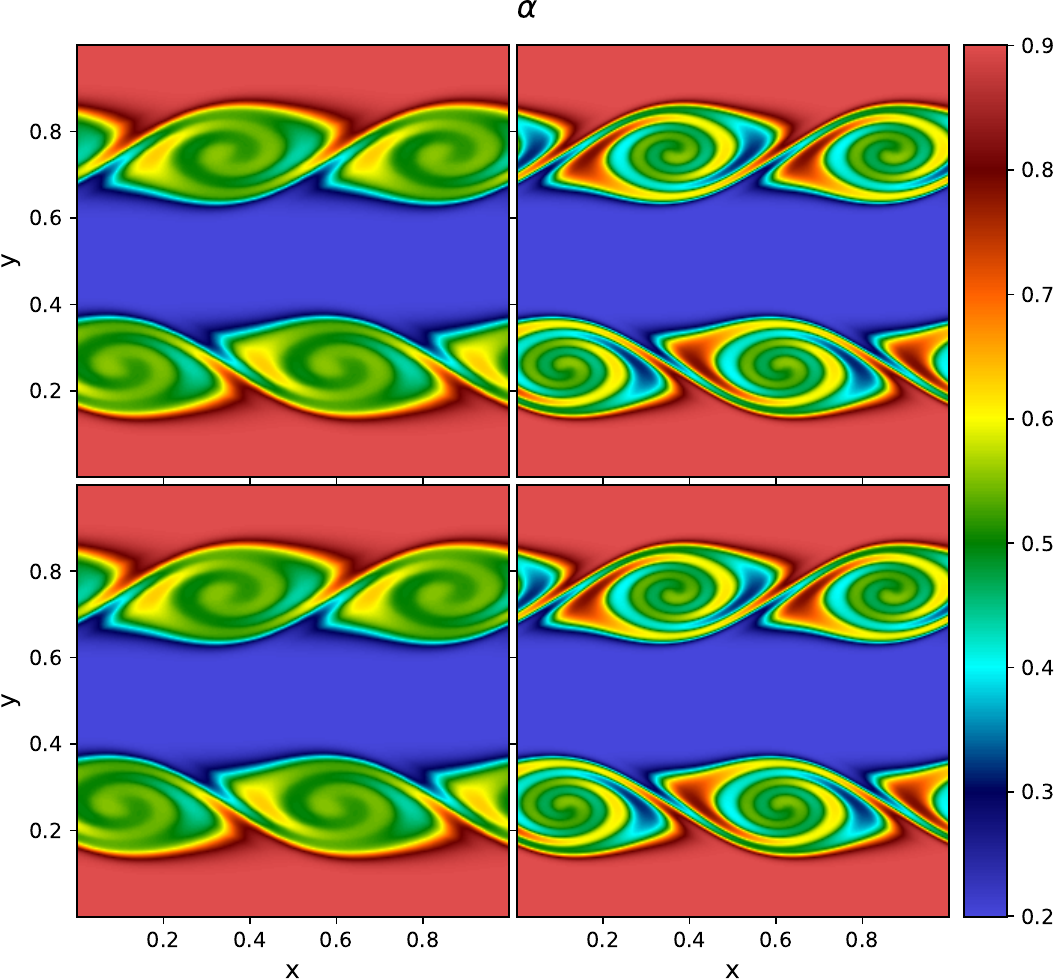}
\caption{Kelvin-Helmholtz instability: numerical solutions for the passively transported volume fraction $\vf_1$ obtained at time $T_f=3$ for the two-fluid single temperature model. The numerical solution is obtained by the new second order RS-IMEX FV scheme. Top panel: $M_{\max} = 10^{-1}$. Bottom panel: $M_{\max} = 3\cdot 10^{-2}$. Left column: $256\times256$ grid. Right column: $512\times512$ grid. }
\label{fig.KH}
\end{figure}

The volume fraction is set as
\begin{equation}
    \vf_1 =
    \begin{cases}
        \vf_L - \vf_m \exp((y-0.25)/L), &\quad \text{if} \quad 0 \leq y < 0.25 \\
        \vf_R + \vf_m \exp(-(y-0.25)/L), &\quad \text{if} \quad 0.25 \leq y < 0.5 \\
        \vf_R + \vf_m \exp((y-0.75)/L), &\quad \text{if} \quad 0.5 \leq y < 0.75 \\
        \vf_L - \vf_m \exp(-(y-0.75)/L), &\quad \text{if} \quad 0.75 \leq y \leq 1 \\
    \end{cases}
\end{equation}
where $\vf_L = 0.9$, $\vf_R = 0.2$ and $\vf_m = (\vf_R - \vf_L)/8$.
In Figure \ref{fig.KH} numerical solutions computed by the second order RS-IMEX FV scheme for the passively transported volume fraction for the Mach numbers $10^{-1}$ and $3 \cdot 10^{-2}$ are depicted. Two different grids consisting of $256\times256$ and $512\times512$ mesh cells and the material CFL condition \eqref{eqn.CFL.cond} $\nu = 0.25$ were used.
The final time is $T_f=3$.
One can observe that despite we deal with a mixture of inviscid fluids, the mesh refinement does
not yield new small scale vortices.
The latter is typical for the Kelvin-Helmholtz instabilities
in an ideal fluid. This is due to the fact that we solve numerically the non-homogeneous system \eqref{eq.SHTC.div.nd}, i.e. physical dissipation is included due to the relative velocity equation. Therefore,
only large vortices are present which corresponds to the frequency modes of the initial data.
Moreover, since the initial data are well-prepared in the sense of \eqref{def.wp}, the $L^1$  errors in the phase densities decreases with the Mach number.
We refer to Table \ref{tab.KH_rho_error} that validates the AP property of the RS-IMEX FV scheme.
\begin{table}[b!]
    \centering
\begin{tabular}{rcc}
    $M\quad$ & $\rho_1$ & $\rho_2$ \\\hline \\[-10pt]
    $10^{-1}$& $1.896\cdot 10^{-3}$ & $1.327\cdot 10^{-3}$ \\
    $3 \cdot 10^{-2}$& $5.037\cdot 10^{-4}$ & $3.525\cdot 10^{-4}$ \\
\end{tabular}
\caption{Kelvin-Helmholtz instability: $L^1$ error of the phase densities for different Mach numbers computed on a mesh with $512\times 512$ grid cells at final time $T_f=3$.}
\label{tab.KH_rho_error}
\end{table}


\section{Conclusions}

We have derived an analyzed a new implicit-explicit finite volume (RS-IMEX FV) scheme for single-temperature SHTC model.
We note that the two-fluid model allows two velocities and pressures.
Further, it includes two dissipative mechanisms: phase pressure and velocity relaxations.
In the proposed scheme these are treated differently.
The relative velocity relaxation
term
is linear and is resolved as a part of the implicit sub-system, whereas the pressure relaxation is
strongly nonlinear and therefore is treated separately by the Newton method.

Our RS-IMEX FV method is constructed in such a way that acoustic-type waves are linearized around a suitably chosen reference state (RS) and approximated implicitly in time and by means of central finite differences in space.
The remaining advective-type waves are approximated explicitly in time and by means of the Rusanov FV method.
The RS-IMEX FV scheme is suitable for all Mach number flows, but in particular it is asymptotic preserving in the low Mach number flow regimes.

Many multi-phase flows, such as granular or sediment transport flows,  can be modeled within the single-temperature approximation.  In turn, many of these flows are weakly compressible and therefore impose severe time step restrictions if solved with a time-explicit numerical scheme.
Therefore, the proposed RS-IMEX FV scheme is suitable to model various environmental flows.

The proposed method was tested on a number of test cases for low and moderately high Mach number flows demonstrating the capability of the scheme to properly capture both regimes.
The theoretical second order accuracy of the scheme was confirmed on a stationary vortex test case.
We compared the second order scheme against its first order variant and shown that the second order scheme has remarkably more accurate approximation of discontinuities.
Finally, the asymptotic preserving property was verified by approximating the Kelvin-Helmholtz instability with well-prepared initial data.

\paragraph{Acknowledgments}
A.T. and M.L. have been partially supported by the Gutenberg Research College, JGU Mainz.
Further, M.L. is grateful for the support of the Mainz Institute of Multiscale Modelling.
I. P. is a member of the Gruppo Nazionale per il Calcolo Scientifico of the Istituto Nazionale di
Alta Matematica (INdAM GNCS) and acknowledges the financial support received from the Italian
Ministry of Education, University and Research (MIUR) in the frame of the Departments of Excellence
Initiative 2018–2022 attributed to the Department of Civil, Environmental and Mechanical
Engineering (DICAM) of the University of Trento (Grant No. L.232/2016).

\noappendicestocpagenum 
\begin{appendices}

\section{Proof of Theorem 4.1}
\label{App:APproof}

We will show the AP property for the first order time semi-discrete scheme \eqref{eq.time-semi.expl} - \eqref{eq.implicit_sub_part}.
Indeed, to obtain a consistent approximation of the limit equations, an appropriate time discretization is essential.
Thereby we will use techniques that were developed in the context of the AP proof for the Euler equations, see for instance \cite{BisLukYel2017}, and for the isentropic two-phase subsystem, see \cite{LukPupTho2022}.
For simplicity we consider without loss of generality $\rrat{1} =1$ and $\rrat{2}=1$.

Let the initial data be well-prepared, i.e. $\q^0 \in \Omega_{wp}$ as given in Definition \ref{def.wp}.
We assume that at time level $t^n$ we have the Mach number expansion for each phase $l = 1,2$
\begin{equation}
    \label{eq.wpdata.tn}
    \rho\indd{l}^n = \rho\indd{l,\RS} + \mathcal{O}(M^2), \quad T^n = T_\RS + \mathcal{O}(M^2) \quad \rho\indd{l,\RS}, T_\RS = \text{ const.},
\end{equation}
in pressure equilibrium up to order $\mathcal{O}(M^3)$, see Assumption \ref{Ann.vf},
\begin{equation}
    \label{eq.wpdata.tn.p}
    p\indd{1,(k)}^n = p\indd{1,(k)}^n, \quad k = 0, 2, \quad p^n\indd{1,(0)} = p\indd{1}(\rho\indd{1,\RS},T_\RS), ~~~p\indd{2,(0)}^n= p\indd{2}(\rho\indd{2,\RS},T_\RS) = \text{ const.}
\end{equation}
Further, for the velocities we assume in concordance with Definition \ref{def.wp} that
\begin{equation}
    \label{eq.wpdata.tn.vel}
    \quad \vv^n = \vv\indu{n}_{(0)} + \mathcal{O}(M^2), \quad \div u\indd{(0)}^n = 0, \quad \ww^n = \mathcal{O}(M^2), \quad \tau^{(\ww)} = M^2.
\end{equation}
Moreover, we assume that the data at the next time level has the Mach number expansion \eqref{eq.Mach.expansion.phase} for the phase densities, temperature and mixture velocity leading to the Mach number expansions \eqref{eq.Mach.expansion.EOS} for pressures, chemical potentials and internal phase energies.
Our aim is to show that the first order IMEX FV method yields a consistent approximation of the incompressible Euler system with variable volume fraction \eqref{eq.SHTC.div.limit}.
To obtain this goal, we show that $\q^{n+1} \in \Omega_{wp}^M$, where the divergence free property of the velocity field is fulfilled up to a $\mathcal{O}(\Delta t)$ term.

Plugging the expansion \eqref{eq.wpdata.tn} at level $t^n$ into the explicit update \eqref{eq.time-semi.expl} we directly have for the volume fraction
\begin{equation}
    \label{AP.vf}
    \vf\indd{1}^{(1)} = \vf\indd{1}^n - \Delta t \vv\indd{(0)}^n \cdot \nabla \vf\indd{1}^n.
\end{equation}
Rewriting equations for $(\vf\indd{1}\rho\indd{1})^{(1)}, (\vf\indd{2}\rho\indd{2})^{(1)}$ in terms of $\rho \vv$ and $\ww$ and using \eqref{AP.vf}, $\div \vv_{(0)}^n = 0$ and $\ww_{(0)}^n = 0$, we have at leading order
    \begin{eqnarray}
        &&  \vf\indd{1}^{(1)}\rho\indd{1,{(0)}}^{(1)} = \vf\indd{1}^n\rho\indd{1,\RS} - \Delta t\rho\indd{1,\RS} \vv\indd{(0)}^n \cdot \nabla \vf\indd{1}^n - \Delta t \div (\rho\indd{(0)}^n\mf\indu{1}\mf\indu{2} \ww^n\indd{(0)}) = \vf\indd{1}^{(1)}\rho\indd{1,\RS},
    \end{eqnarray}
thus $\rho\indd{1,(0)}\indu{(1)} = \rho\indd{1,\RS}$.
With the same strategy we obtain $\rho\indu{n+1}\indd{1,(1)} = 0$.
Analogously, we obtain $\rho\indd{2,(0)}\indu{(1)} = \rho\indd{2,\RS}$ and $\rho\indu{(1)}\indd{2,(1)} = 0$.
Summarizing, the phase densities satisfy the expansion \eqref{eq.wpdata.tn} at the intermediate time level $t^{(1)}$.
Using $\ww_{(0)}^n = 0$ and the evolution of the volume fraction \eqref{AP.vf}, we obtain for the momentum and relative velocity equations
\begin{subequations}
    \label{AP.v_w_explicit}
    \begin{eqnarray}
        &&(\rho\vv)_{(0)}^{(1)} = (\rho\vv)_{(0)}^n - \Delta t \vv\indd{(0)}^n \cdot \nabla (\rho\vv)\indd{(0)}^n, \\
        && \ww_{(0)}^{(1)} = 0, \quad \ww_{(1)}^{(1)} = 0.
    \end{eqnarray}
\end{subequations}
Multiplying the energy equation in the explicit update \eqref{eq.time-semi.expl} by $M^2$ and using the notation \eqref{eq.scaled_total_energy}, yields
\begin{equation}
     (\vf\indd{1}\rho\indd{1}e\indd{1})^{(1)} + (\vf\indd{2}\rho\indd{2}e\indd{2})^{(1)} + M^2 (\rho E\indd{\tkin})^{(1)} = (\vf\indd{1}\rho\indd{1}e\indd{1})^{n} + (\vf\indd{2}\rho\indd{2}e\indd{2})^{n} + M^2 (\rho E\indd{\tkin})^{n} + \mathcal{O}(M^2).
\end{equation}
For the leading order terms of the internal energy, we obtain directly $$\left(\vf\indd{1}(\rho\indd{1}e\indd{1})\indd{(0)} + \vf\indd{2}(\rho\indd{2}e\indd{2})\indd{(0)}\right)^{(1)} = \left(\vf\indd{1}(\rho\indd{1}e\indd{1})\indd{(0)} + (\vf\indd{2}(\rho\indd{2}e\indd{2})\indd{(0)}\right)^n$$ which completes the analysis of the explicit part \eqref{eq.time-semi.expl}.

For the implicit part, we will follow the reasoning of \cite{BisLukYel2017}, where the AP property is shown for the Euler equations analyzing the structure of the implicit elliptic operator.
Since $\alpha\indd{1}\rho\indd{1}$ does not change during the implicit part, the expansion of the phase densities at $t^{n+1}$ fulfills \eqref{eq.wpdata.tn}.
Therefore, we obtain $\rho\indd{1,(0)}\indu{n+1} = \rho\indd{1,\RS}$ and analogously $\rho\indd{2,(0)}\indu{n+1} = \rho\indd{2,\RS}$ for the second phase density.

Next, we analyze the elliptic update of the total energy \eqref{eq.impl.Etot}.
Analogously to the fully discrete operators $\mathcal{L}_I$ and $\mathcal{K}_I$ in \eqref{eq.Elliptic_fully_disc}, we define semi-discrete operators
\begin{eqnarray}
    &&L_h = \div \left(\left(\frac{(\vf\indd{1}\rho\indd{1}e\indd{1})^n + (\vf\indd{2}\rho\indd{2}e\indd{2})^{n} + M^2 (\rho E\indd{\tkin})^{n} + \vf_1^n p_1^n + \vf_2^n p_2^n}{\rho^{n+1}}\right)\nabla (\phi_p^n - 1)\right),\\
    &&K_h = \div \left(\tau^{(\ww)} \frac{(\mu_1^n - \mu_2^2)(\rho \mf_1 \mf_2)^{n+1}}{\tau^{(\ww)} + \Delta t (\mf_1 \mf_2)^{n+1}}\nabla \left(\partial \mu_\RS^n\right)\right).
\end{eqnarray}
Note that with \eqref{eq.p.expandsion.E} we have $L_h = \mathcal{O}(1).$ From \eqref{eq.dmu.exp} and $\tau^{(\ww)} = \mathcal{O}(M^2)$ it follows $K_h = \mathcal{O}(M^2)$.
Using the notation as in \eqref{eq.scaled_total_energy}, we define
\begin{equation}
    \rho \check{E} = \left(\vf_1 \rho_1 e_1 + \vf_2 \rho_2 e_2 + M^2 \rho E_\tkin\right),\quad \check{p} = (\vf_1 p_1 + \vf_2 p_2)\quad \check{\mu} = \mu_1 - \mu_2.
\end{equation}
Now taking into account the scaling of $\ww^n$ given in \eqref{eq.wpdata.tn.vel}, we write the implicit update for the total energy \eqref{eq.Elliptic_fully_disc} as
\begin{equation}
    \label{eq.Semi-discrete-operators}
    \left(\bm{I} - \frac{\Delta t^2}{M^2} (L_h + K_h)\right) (\rho \check{E})^{n+1} = (\rho \check{E})^n - \Delta t\nabla_h \cdot \left((\rho \check{E} + \check{p})^n \vv^{(1)}\right) - \frac{\Delta t^2}{M^2} L_h (M^2\rho E_\tkin^n) - \frac{\Delta t^2}{M^2} K_h (\rho \check{E})_\RS^n + \mathcal{O}(M^2).
\end{equation}
The operators $L_h$ and $K_h$ are symmetric, positive definite and the inverse  of $A = \bm{I} - \frac{\Delta t^2}{M^2} (L_h + K_h)$ exists.
Consequently, system \eqref{eq.Semi-discrete-operators} has a unique solution for any $M>0$.
Similar as in \cite{BisLukYel2017}, we obtain that the eigenvalues of $A^{-1}$ are $1$ and $\mathcal{O}(M^2)$. Applying analogous arguments as in \cite[Lem. 4.6]{BisLukYel2017}, we derive
\begin{equation}
    (\rho e)^{n+1} = (\rho e)^n - \Delta t\nabla_h \cdot \left(\left(\rho e + \alpha_1 p_1 + \alpha_2 p_2\right)^n \vv^{n}\right) + \mathcal{O}(M^2).
\end{equation}
Focusing on the leading order terms and using the evolution of the volume fraction \eqref{AP.vf}, $\div \vv_{(0)}^{n} = 0$, see \eqref{eq.wpdata.tn.vel}, $p_{1,(0)}^n = p_{2,(0)}^n$, see \eqref{eq.wpdata.tn.p}, and EOS \eqref{eq.EOS},  yields for the temperature
the following expansion
\begin{align*}
    \begin{split}
        \left(\vf_1^{(1)}\rho_{1,\RS}c_{v,1} + \vf_2^{(1)}\rho_{2,\RS} c_{v,2}\right)T^{n+1}_{(0)}=& \left(\vf_1^{n}\rho_{1,\RS}c_{v,1} + \vf_2^{n}\rho_{2,\RS} c_{v,2}\right)T_{\RS} - \Delta t \vv_{(0)}^n \left(\rho_{1,\RS}c_{v,1} - \rho_{2,\RS} c_{v,2}\right) T_\RS \nabla \vf_{1}^n + \mathcal{O}(M^2)\\
        =&  \left(\vf_1^{(1)}\rho_{1,\RS}c_{v,1} + \vf_2^{(1)}\rho_{2,\RS} c_{v,2}\right)T_\RS + \mathcal{O}(M^2).
    \end{split}
\end{align*}
Since the factor $\vf_1^{(1)}\rho_{1,\RS}c_{v,1} + \vf_2^{(1)}\rho_{2,\RS} c_{v,2}$ is positive and independent of $M$, we derive $T_{(0)}^{n+1} = T_\RS + \mathcal{O}(M^2)$, thus the temperature has a correct asymptotic expansion.
Moreover, in the limit as $M \to 0$, we obtain $T = T_\RS$.
Further, the update of the relative velocity \eqref{eq.impl.w} and momentum \eqref{eq.impl.mom} yield
\begin{align}
    \ww_{(0)}^{n+1} &= 0, \quad \nabla \mu_{(0)}^{n+1} = 0, \quad \nabla \mu_{(1)}^{n+1} = 0,\\
    (\rho\vv)_{(0)}^{n+1} &= (\rho\vv)_{(0)}^{(1)} - \Delta t \nabla p_{(2)}^{n+1}, \quad \nabla p_{(0)}^{n+1} = 0, \quad \nabla p_{(1)}^{n+1} = 0.
\end{align}
Since the mass densities of the phases are not evolved at the implicit step, it holds $\rho^{n+1} = (\vf_1\rho_1)^{(1)} + (\vf_2\rho_2)^{(1)}$.
Using the volume fraction, we can rewrite the momentum equation as
\begin{equation}
    \vv_{(0)}^{n+1} = \vv_{(0)}^{(1)} - \Delta t \frac{\rho^n}{\rho^{n+1}}\vv_{(0)}^n \cdot \nabla \vv_{(0)}^n - \Delta t \frac{\nabla p_{(2)}^{n+1}}{\rho^{n+1}}
\end{equation}
which is consistent with the low Mach number limit \eqref{eq.SHTC.div.limit} up to a $\mathcal{O}(\Delta t)$ term.
From the energy equation \eqref{eq.implicit_sub_part_energy} we obtain
\begin{equation}
    \label{eq.AP_proof_Mach_exp_int_e_new}
    (\rho e)_{(0)}^{n+1} = (\rho e)^n_{(0)} - \Delta t \div\left(\left((\vf_1^n (\rho e)_{1,\RS} + \vf_2^n (\rho e)_{2,\RS}) + p_{(0)}^n\right) \vv^{n+1}_{(0)}\right) + \mathcal{O}(M^2).
\end{equation}
Using the definition of the internal mixture energy, we obtain
\begin{align*}
    \vf_1^{(1)} (\rho e)_{1,\RS} + \vf_2^{(1)} (\rho e)_{2,\RS} =&~ \vf_1^n (\rho e)_{1,\RS} + \vf_2^n (\rho e)_{2,\RS} \\
    &~- \Delta t \div\left(\left((\vf_1^n (\rho e)_{1,\RS} + \vf_2^n (\rho e)_{2,\RS}) + p_{(0)}^n\right) \vv^{n+1}_{(0)}\right) + \mathcal{O}(M^2).
\end{align*}
Applying the evolution of the volume fraction \eqref{AP.vf} and $\nabla p_{(0)}^n = 0$, we obtain
\begin{align*}
    (\vf\indd{1}^n - \Delta t \vv\indd{(0)}^n \cdot \nabla \vf\indd{1}^n)(\rho e)_{1,\RS} + (\vf_2^n + \Delta t \vv\indd{(0)}^n \cdot \nabla \vf\indd{1}^n) (\rho e)_{2,\RS} =&~ \vf_1^n (\rho e)_{1,\RS} + \vf_2^n (\rho e)_{2,\RS} \\
    &~- \Delta t \vv^{n+1}_{(0)} \cdot \nabla \left((\vf_1^n (\rho e)_{1,\RS} + \vf_2^n (\rho e)_{2,\RS})\right)\\
    &~- \Delta t \left((\vf_1^n (\rho e)_{1,\RS} + \vf_2^n (\rho e)_{2,\RS}) + p_{(0)}^n\right) \div \vv^{n+1}_{(0)} + \mathcal{O}(M^2),
\end{align*}
which reduces to
\begin{align*}
    (\vv^{n}_{(0)} - \vv^{n+1}_{(0)}) \cdot \nabla \left(\vf_1^n (\rho e)_{1,\RS} + \vf_2^n (\rho e)_{2,\RS}\right) = \left((\vf_1^n (\rho e)_{1,\RS} + \vf_2^n (\rho e)_{2,\RS}) + p_{(0)}^n\right) \div \vv^{n+1}_{(0)} + \mathcal{O}(M^2).
\end{align*}
The left hand side is of order $\mathcal{O}(\Delta t)$ and the factor $(\vf_1^n (\rho e)_{1,\RS} + \vf_2^n (\rho e)_{2,\RS}) + p_{(0)}^n$ is positive. Consequently, we obtain the result $\div \vv_{(0)}^{n+1} = \mathcal{O}(\Delta t) + \mathcal{O}(M^2)$.

Finally, we apply the pressure relaxation on the volume fraction, and we obtain
\begin{equation}
    p_1^{n+1} = p_2^{n+1} - \tau^{(\vf)} \rho^{n+1} M^2 \left(\frac{\vf^{n+1}- \vf^{n}}{\Delta t}\right).
\end{equation}
Thus $p_{1,{(0)}}^{n+1} = p_{2,{(0)}}^{n+1}$ and $ p_{1,{(2)}}^{n+1} = p_{2,{(2)}}^{n+1}$ since $p_{1,{(1)}}^{n+1} = p_{2,(1)}^{n+1} = 0$.
Note that this is due to phase densities and the temperature fulfill expansion \eqref{eq.wpdata.tn} at the new time level $t^{n+1}$.
Moreover, $\rho^{n+1}(\frac{\vf^{n+1}- \vf^{n}}{\Delta t}) = \mathcal{O}(1)$ with respect to the Mach number, since the time step is independent of the Mach number as well.

\noindent This concludes the proof.

	\section{Polar coordinates}
\label{App:Polar}
We consider a continuous solution of the homogeneous part of system \eqref{eq.SHTC.div} without relaxation source terms.
Let the Cartesian coordinates in 2D be denoted by $x = (x_1,x_2)$.
We define the polar coordinates in terms of radius $r$ and angle $\theta$ as
\begin{equation}
    x_1 = r \cos(\theta), \quad x_2 = r \sin(\theta).
\end{equation}
The velocity based quantities are defined by
\begin{subequations}
    \begin{eqnarray}
        && v_1 = v_r \cos(\theta) - v_\theta \sin(\theta), \quad w_1 = w_r \cos(\theta) - w_\theta \sin(\theta),\\
        && v_2 = v_r \cos(\theta) + v_\theta \sin(\theta), \quad w_2 = w_r \cos(\theta) + w_\theta \sin(\theta).
    \end{eqnarray}
\end{subequations}
Using
\begin{equation}
    \frac{\partial x_1}{\partial r} = \cos(\theta), \quad \frac{\partial x_1}{\partial \theta} = -\frac{\sin(\theta)}{r}, \quad \frac{\partial x_2}{\partial r} = \sin(\theta), \quad \frac{\partial x_1}{\partial \theta} = \frac{\cos(\theta)}{r},
\end{equation}
we obtain for
\begin{equation}
    \q = (\alpha^1, \alpha^1 \rho^1, \alpha^2 \rho^2, v_r, v_\theta, w_r, w_\theta, \rho E)^T
\end{equation}
the following system in polar coordinates
\begin{subequations}
    \label{eq:TFM_polar}
    \begin{eqnarray}
        && \frac{\partial \alpha^1}{\partial t}+\frac{\ur}{r}\frac{\partial}{\partial r} \left(r \alpha^1 \right) + \frac{\uphi}{r}\frac{\partial}{\partial \theta} \alpha^1 = 0, \\
        && \frac{\pd (\alpha^1\rho^1)}{\pd t}+\frac{1}{r}\frac{\pd}{\pd r} \left(r \alpha^1 \rho^1 \ur^1 \right) + \frac{1}{r}\frac{\pd}{\pd \theta} \left( \vf\rho^1 \uphi^1 \right) =0,\\[2mm]
        && \frac{\pd (\alpha^2\rho^2)}{\pd t}+\frac{1}{r}\frac{\pd}{\pd r} \left(r \alpha^2 \rho^2 \ur^2 \right) + \frac{1}{r}\frac{\pd}{\pd \theta} \left(  \alpha^2 \rho^2 \uphi^2 \right) =0,\\[2mm]
        &&
        \begin{split}
            \frac{\pd (\rho \ur)}{\pd t}& +\frac{1}{r}\frac{\pd}{\pd r}\left(r\left(\rho \ur^2 + \rho \mf (1-\mf) w_r w_r + p\right)\right) \\
            & + \frac{1}{r}\frac{\pd}{\pd \theta}\left(\rho \ur \uphi +  \rho \mf c^2 w_r \wphi\right) = \frac{\rho \uphi^2 + \rho \mf (1-\mf) \wphi^2 + p}{r},
        \end{split}\\
        &&
        \begin{split}
            \frac{\pd (\rho \uphi)}{\pd t}& +\frac{1}{r}\frac{\pd}{\pd r}\left(r\left(\rho \ur \uphi + \rho \mf (1-\mf) w_r \wphi\right)\right) \\
            & + \frac{1}{r}\frac{\pd}{\pd \theta}\left( \rho \uphi^2 + p + \rho \mf (1-\mf) \wphi^2\right) = -\frac{\rho \ur \uphi + \rho \mf c^2 w_r \wphi }{r},
        \end{split}\\
        &&
        \begin{split}
            \frac{\pd w_r}{\pd t}&+\frac{\pd}{\pd r} \left(\ur w_r + \uphi \wphi + (1-2\mf) \frac{w_r w_r + \wphi\wphi}{2} + \mu_1 - \mu_2\right)\\
            & + \uphi \left(\frac{1}{r}\frac{\pd}{\pd \theta} w_r - \frac{1}{r} \frac{\pd}{\pd r}\left(r w_\theta\right)\right) = 0,
        \end{split}\\[2mm]
        &&
        \begin{split}
            \frac{\pd \wphi }{\pd t}&+\frac{1}{r}\frac{\pd}{\pd \theta} \left(\ur \wr + \uphi \wphi + (1-2\mf) \frac{\wr \wr + \wphi \wphi}{2} + \mu_1 - \mu_2 \right)\\
            &+\ur \left(\frac{1}{r} \frac{\pd}{\pd r}\left(r w_\theta\right) - \frac{1}{r}\frac{\pd}{\pd \theta} w_r\right)  = 0,
        \end{split}\\[2mm]
        &&
        \begin{split}
            \frac{\pd (\rho E)}{\pd t}& + \frac{1}{r}\frac{\pd}{\pd r} \left( r\left(\ur (\rho E + p)  +\rho ~\left[\ur \wr + \uphi \wphi +
            (1-2\mf\indu{1}) \frac{\wr^2 + \wphi^2}{2}\right]\mf^1\mf^2 \wr \right)\right) \\
            & + \frac{1}{r}\frac{\pd}{\pd \theta}\left(\uphi (\rho E + p)  +\rho ~\left[\ur \wr + \uphi \wphi +
            (1-2\mf\indu{1}) \frac{\wr^2 + \wphi^2}{2}\right]\mf^1\mf^2 \wphi \right) =0.
        \end{split}
    \end{eqnarray}
\end{subequations}


\end{appendices}

\bibliographystyle{plain}
\bibliography{lit_twofluid.bib}
%

\end{document}